\documentclass{amsart}

\usepackage[utf8]{inputenc}
\usepackage[T1]{fontenc}
\usepackage{ae,aecompl}
\usepackage{subfigure}

\usepackage{amssymb,amscd}
\usepackage{algorithmic}
\usepackage[matrix,arrow,curve]{xy}

\usepackage{tikz}
\usetikzlibrary{arrows}

\usepackage{verbatim}

\title[Hecke group algebras and their representation theory]
{The Hecke group algebra of a Coxeter group\\ and its representation theory}

\author{Florent Hivert and Nicolas M.~Thiéry}

\address{Florent Hivert -- \rm\texttt{florent.hivert@univ-rouen.fr}
 -- LIFAR -- university of Rouen --  Avenue de l'Université --
  Technopôle du Madrillet --
76801 Saint Etienne du Rouvray cedex -- FRANCE}

\address{Nicolas M. Thiéry -- \rm\texttt{nthiery@users.sf.net}
  Univ Paris-Sud, Laboratoire de Mathématiques d'Orsay,
  Orsay, F-91405; CNRS, Orsay, F-91405}



\newcommand{\TODO}[1]{}
\newcommand{\FIXME}[1]{}

\renewcommand{\k}{{\mathbb{C}}} 
\newcommand{\NN}{{\mathbb{N}}}
\newcommand{\F}{{\mathbb{F}}}
\newcommand{\CC}{{\mathbb{C}}}

\newcommand{\s}{\sigma}
\renewcommand{\t}{\tau}
\newcommand{\End}{{\operatorname{End}}}
\newcommand{\sym}{{\operatorname{Sym}}}
\newcommand{\qsym}{{\operatorname{QSym}}}
\newcommand{\ncsf}{{\operatorname{NCSF}}}

\newcommand{\W}{{W}}
\newcommand{\Wa}[1][\CC]{{#1[\W]}}
\newcommand{\kW}[1][\CC]{{#1\W}}
\newcommand{\heckeW}[2][\W]{{\mathcal{H}(#1)(#2)}}
\newcommand{\affineheckeW}[2][\W]{{\widetilde{\operatorname{H}}(#1)(#2)}}
\newcommand{\heckeWW}[1][\W]{{\mathcal{H}#1}}

\newcommand{\sg}[1]{{\mathfrak{S}_{#1}}}
\newcommand{\sga}[2][\CC]{{#1[\sg{#2}]}}
\newcommand{\ksg}[2][\CC]{{#1\sg{#2}}}
\newcommand{\hecke}[2]{{\mathcal{H}_{#1}(#2)}}
\newcommand{\heckesg}[1]{\heckeWW[\sg{#1}]}
\newcommand{\affinehecke}[2][n]     {{\widetilde{\operatorname{H}}_{#1}(#2)}}
\newcommand{\TL}[2]{\operatorname{\mathcal{T\!L}}_{#1}(#2)}

\newcommand{\ndf}[1]{\operatorname{NDF}_{#1}}
\newcommand{\ndfa}[2][\CC]{{#1[\ndf{#2}]}}
\newcommand{\pf}[1]{\operatorname{PF}_{#1}}
\newcommand{\ndpf}[1]{\operatorname{NDPF}_{#1}}
\newcommand{\ndpfa}[2][\CC]{{#1[\ndpf{#2}]}}
\newcommand\gr[1]{\operatorname{gr_{#1}}} 

\newcommand{\initial}[1]{\operatorname{Init}_{#1}}
\newcommand{\ndinitial}[1]{\operatorname{NDInit}_{#1}}

\newcommand{\Des}{{\operatorname{Des}}}
\newcommand{\Rec}{{\operatorname{iDes}}}
\newcommand{\cset}{{\operatorname{C}}}
\newcommand{\partof}{\vdash}                    
\newcommand{\compof}{\vDash}                    
\newcommand{\fin}{{\succeq}}                    
\newcommand{\lon}{\operatorname{\ell}}          
\newcommand{\act}{\cdot}

\newcommand{\melem} {{e}}                       
\newcommand{\ls}[1][s] {{\overleftarrow{#1}}}
\newcommand{\opi}{{\overline \pi}}

\newcommand{\std}{{\operatorname{std}}}
\newcommand{\init}{{\operatorname{init}}}
\newcommand{\id}{{\operatorname{id}}}
\newcommand{\im}{{\operatorname{im}}}
\newcommand{\fibers}{{\operatorname{fibers}}}

\newcommand{\rad}{{\operatorname{rad}}}
\newcommand{\Hom}{{\operatorname{Hom}}}
\newcommand{\cartan}{C}

\newcommand{\calG}{{\mathcal G}}
\newcommand{\calK}{{\mathcal K}}
\newcommand{\suchthat}{{\ |\ }}

\newtheorem{theorem}{Theorem}[section]
\newtheorem{lemma}[theorem]{Lemma}
\newtheorem{proposition}[theorem]{Proposition}
\newtheorem{corollary}[theorem]{Corollary}
\newtheorem{corollaries}[theorem]{Corollaries}
\theoremstyle{definition}
\newtheorem{definition}[theorem]{Definition}
\newtheorem{example}[theorem]{Example}

\newtheorem{conjecture}[theorem]{Conjecture}

\theoremstyle{remark}
\newtheorem{remark}[theorem]{Remark}

\begin{document}
\begin{abstract}
  Let $W$ be a finite Coxeter group. We define its Hecke-group algebra
  by gluing together appropriately its group algebra and its $0$-Hecke
  algebra. We describe in detail this algebra (dimension, several
  bases, conjectural presentation, combinatorial construction of
  simple and indecomposable projective modules, Cartan map) and give
  several alternative equivalent definitions (as symmetry preserving
  operator algebra, as poset algebra, as commutant algebra, ...).


  In type $A$, the Hecke-group algebra can be described as the algebra
  generated simultaneously by the elementary transpositions and the
  elementary sorting operators acting on permutations. It turns out to
  be closely related to the monoid algebras of respectively
  nondecreasing functions and nondecreasing parking functions, the
  representation theory of which we describe as well.

  This defines three towers of algebras, and we give explicitly the
  Grothendieck algebras and coalgebras given respectively by their
  induction products and their restriction coproducts. This yields some
  new interpretations of the classical bases of quasi-symmetric and
  noncommutative symmetric functions as well as some new bases.

\end{abstract}

\maketitle

\tableofcontents

\section{Introduction}

Given an \emph{inductive tower of algebras}, that is a sequence of algebras
\begin{equation}
A_0 \hookrightarrow A_1 \hookrightarrow \cdots \hookrightarrow A_n
\hookrightarrow \cdots,
\end{equation}
with embeddings $A_m\otimes A_n \hookrightarrow A_{m+n}$ satisfying an
appropriate associativity condition, one can introduce two \emph{Grothendieck
rings}
\begin{equation}
\calG(A):=\bigoplus_{n\ge 0}G_0(A_n)\qquad \text{and} \qquad
\calK(A):=\bigoplus_{n\ge 0}K_0(A_n)\,,
\end{equation}
where $G_0(A)$ and $K_0(A)$ are the (complexified) Grothendieck groups
of the categories of finite-dimensional $A$-modules and projective
$A$-modules respectively, with the multiplication of the classes of an
$A_m$-module $M$ and an $A_n$-module $N$ defined by the induction product
\begin{equation}
[M] \cdot [N] = [M\widehat{\otimes} N] =
[M\otimes N \uparrow_{A_m\otimes A_n}^{A_{m+n}}]\,.
\end{equation}

If $A_{m+n}$ is a projective $A_m\otimes A_n$-module, one can define a
coproduct on these rings by means of restriction of representations,
turning these into coalgebras. Under favorable circumstances the
product and the coproduct are compatible turning these into mutually
dual Hopf algebras.

The basic example of this situation is the character ring of the
symmetric groups (over $\CC$), due to Frobenius. Here the
$A_n:=\sga{n}$ are semi-simple algebras, so that
\begin{equation}
{G}_0(A_n) = {K}_0(A_n)= R(A_n)\,,
\end{equation}
where $R(A_n)$ denotes the vector space spanned by isomorphism classes
of indecomposable modules which, in this case, are all simple and
projective.
The irreducible representations $[\lambda]$ of $A_n$ are parametrized by
partitions $\lambda$ of $n$, and the Grothendieck ring is isomorphic to the
algebra $\sym$ of symmetric functions under the
correspondence
$[\lambda] \leftrightarrow s_\lambda$,
where $s_\lambda$ denotes the Schur function associated with
$\lambda$.  Other known examples with towers of group algebras over
the complex numbers $A_n:=\CC[G_n]$ include the cases of wreath
products $G_n := \Gamma\wr\sg{n}$ (Specht), finite linear groups $G_n
:= GL(n,\F_q)$ (Green), \emph{etc.}, all related to symmetric
functions (see~\cite{MacDonald.SF.1995,Zelevinsky.1981}).

Examples involving non-semisimple specializations of Hecke algebras have also
been worked out.
Finite Hecke algebras of type $A$ at roots of unity ($A_n=\hecke{n}{\zeta}$,
$\zeta^r=1$) yield quotients and subalgebras of $\sym$~\cite{Lascoux_Leclerc_Thibon.1996}.
The Ariki-Koike algebras at roots of unity give rise to level $r$ Fock spaces
of affine Lie algebras of type $A$~\cite{Ariki_Koike.1994}.
The $0$-Hecke algebras $A_n=\hecke{n}{0}$ corresponds to the pair
Quasi-symmetric functions / Noncommutative symmetric functions,
${\mathcal G}=\qsym$, ${\mathcal
  K}=\ncsf$~\cite{Krob_Thibon.NCSF4.1997}.  Affine Hecke algebras at
roots of unity lead to the positive halves of the enveloping algebras
$U(\widehat{sl}_r)$ and $U(\widehat{sl}_r)^*$~\cite{Ariki.1996}, and
the case of generic affine Hecke algebras can be reduced to a
subcategory admitting as Grothendieck rings the positive halves of the
enveloping algebras $U(\widehat{gl}_\infty)$ and
$U(\widehat{gl}_\infty)^*$~\cite{Ariki.1996}.
Further interesting examples are the tower of $0$-Hecke-Clifford
algebras~\cite{Olshanski.1992,Bergeron_Hivert_Thibon.2004} giving rise to the peak algebras~\cite{Stembridge.1997},
and a degenerated version of the Ariki-Koike algebras~\cite{HNTAriki}
giving rise to a colored version of $\qsym$ and $\ncsf$.

\bigskip

The original aim of this research was the study of the representation
theories of several towers of algebras related to the symmetric groups
and their Hecke algebras $\hecke{n}{q}$, in order to derive new
examples of Grothendieck algebras and coalgebras. Along the way, one
of these algebras became interesting in itself and for all types. Here
is the structure of the paper together with the main results.

In Section~\ref{section.heckesg}, we introduce the main object of this
paper, namely the \emph{Hecke-Group algebra} $\heckeWW$ of a (finite)
Coxeter group $\W$. It is constructed as the smallest algebra
containing simultaneously the group algebra of $\W$ and its $0$-Hecke
algebra (and in fact any other $q$-Hecke algebra of $\W$). It turns
out that this algebra has unexpectedly nice properties. We first show
that $\heckeWW$ is better understood as the algebra of antisymmetry
(or symmetry) preserving operators; this allows us to compute its
dimension, and to give explicit bases. We further realize it as the
incidence algebra of a pre-order and derive from it its representation
theory. In particular, we construct explicitly the projective and
simple modules. The Cartan matrix suggests a link between $\heckeWW$
and the incidence algebra of the boolean lattice.  We actually show
that these algebras are Morita equivalent.

Turning back to type $A$, we get a new tower of algebras
$\heckesg{n}$. Specifically, each $\heckesg{n}$ is the algebra
generated by both elementary transpositions and elementary sorting
operators acting on permutations of $\{1,\dots,n\}$. We compute the
restrictions and inductions of simple and projective modules. This
gives rise to a new interpretation of some bases of quasi-symmetric
and noncommutative symmetric functions in representation theory.

In Sections~\ref{section.ndf} and ~\ref{section.ndpf} we turn to the
study of two other towers, namely the towers of the monoid algebras of
nondecreasing functions $\ndfa n$ and of nondecreasing parking
functions $\ndpfa n$.  We prove that those two algebras are the
respective quotients of $\heckesg{n}$ and $\hecke n{0}$, through their
representations on exterior powers. We deduce the structure of their
projective and simple modules, their Cartan matrices, and the
induction and restriction rules.  We also show that the algebra of
nondecreasing parking functions is isomorphic to the incidence algebra
of some lattice.

The following diagram summarizes the relations between all the
aforementioned towers of algebras, as well as the Temperley-Lieb
algebra $\TL{n}{-1}$ at $q=-1$ and the respective quotients $\sga{n}
\circlearrowleft\!  \bigwedge\!\k^n$ and
$\hecke{n}{q}\circlearrowleft\!  \bigwedge\!\k^n$ of $\sga{n}$ and
$\hecke{n}{q}$ through their representations on exterior powers:
\begin{equation}
  \vcenter{
  \xymatrix@R=1cm@C=0.50cm{
    \hecke{n}{-1}        \ar@{->>}[d] \ar@{^(->}@/^4ex/[rrrr] &
    \hecke{n}{0}         \ar@{->>}[d] \ar@{^(->}@/^3ex/[rrr]  &
    \hecke{n}{1}=\sga{n} \ar@{->>}[d] \ar@{^(->}@/^2ex/[rr]   &
    \hecke{n}{q}         \ar@{->>}[d] \ar@{^(->}    [r]       &
    H\sg{n}              \ar@{->>}[d]\\
    \TL{n}{-1}           \ar@{^(->}@/_4ex/[rrrr]&
    \ndpfa n                                         \ar@{^(->}@/_3ex/[rrr] &
    \sga{n}     \circlearrowleft\! \bigwedge\!\k^n \ar@{^(->}@/_2ex/[rr]  &
    \hecke{n}{q}\circlearrowleft\! \bigwedge\!\k^n \ar@{^(->}    [r]      &
    \ndfa n
  }}
\end{equation}
\smallskip

Finally, in Section~\ref{section.researchInProgress}, we discuss
further research in progress, in particular about the links between
$\heckeWW$ and the affine Hecke algebras.

\subsubsection*{Acknowledgments}

We would like to thank Jean-Yves Thibon who is at the origin of the
present research, and for numerous fruitful discussions.

This paper mostly reports on computation driven research using the
package \texttt{MuPAD-Combinat} by the authors of the present
paper~\cite{MuPAD-Combinat}. This package is designed for the computer
algebra system \texttt{MuPAD} and is freely available from
\texttt{http://mupad-combinat.sf.net/}. Among other things, it allows
to automatically compute the dimensions of simple and indecomposable
projective modules together with the matrix of Cartan invariants of a
finite dimensional algebra, knowing its multiplication table. We also
would like to thank the \texttt{MuPAD} group for their support in the
development of \texttt{MuPAD-Combinat}.

\section{Background}

\subsection{Compositions and sets}

Let $n$ be a fixed integer. Recall that each subset $S$ of
$\{1,\dots,n-1\}$ can be uniquely identified with a $p$-tuple $I :=
(i_1,\dots,i_p)$ of positive integers of sum $n$:
\begin{equation}
  \label{eq.set.to.comp}
  S=\{s_1 < s_2 < \dots < s_p\}
  \longmapsto \cset(S):=(s_1,s_2-s_1,s_3-s_2,\ldots,n-s_p)\,.
\end{equation}
We say that $I$ is a \emph{composition of $n$} and we write it by
$I\compof n$. The converse bijection, sending a composition to its
\emph{descent set}, is given by:
\begin{equation}
  \label{eq.comp.to.set}
  I = (i_1,\dots,i_p) \longmapsto
  \Des(I) = \{i_1+\dots+i_j \suchthat j=1,\dots,p-1\}\,.
\end{equation}
The number $p$ is called the \emph{length} of $I$ and is denoted by $\lon(I)$.


The notions of complement of a set $S^c$ and of inclusion of sets
can be transfered to compositions, leading to the complement of a
composition $K^c$ and to the refinement order on compositions: we say
that $I$ is {\em finer} than $J$, and write $I\fin J$, if and only if
$\Des(I)\supseteq \Des(J)$.


\subsection{Coxeter groups and Iwahori-Hecke algebras}

Let $(\W, S)$ be a Coxeter group, that is a group $W$ with a
presentation
\begin{equation}
  W = \langle\, S\ \suchthat\  (ss')^{m(s,s')},\ \forall s,s'\in S\,\rangle\,,
\end{equation}
where each $m(s,s')$ is in $\{1,2,\dots,\infty\}$, and $m(s,s)=1$. 
The elements $s\in S$ are called \emph{simple reflections}, and the
relations can be rewritten as:
\begin{equation}
  \begin{alignedat}{2}
    s^2 &=\id\,, &\quad& \text{ for all $s\in S$}\,,\\
    \underbrace{ss'ss's\cdots}_{m(s,s')} &=
    \underbrace{s'ss'ss'\cdots}_{m(s,s')}\,, && \text{ for all $s,s'\in S$}\ ,
  \end{alignedat}
\end{equation}
Most of the time, we just write $W$ for $(W,S)$. In general, we follow
the notations from~\cite{Bjorner_Brenti.2005}, and we refer to this
monograph for details on Coxeter groups and their Hecke algebras.
Unless stated otherwise, we always assume that $\W$ is finite.

The prototypical example of Coxeter group is the $n$-th symmetric
group $(W,S) := A_{n-1} = (\sg{n}, \{s_1,\dots,s_{n-1}\})$, where
$s_i$ denotes the \emph{elementary transposition} which exchanges $i$
and $i+1$. The relations are given by:
\begin{equation}
  \begin{alignedat}{2}
    s_i^2           & = \id\,,                &     & \text{ for } 1\leq i\leq n-1\,,\\
    s_i s_j         & = s_j s_i\,,            &     & \text{ for } |i-j|\geq2\,, \\
    s_i s_{i+1} s_i & = s_{i+1} s_i s_{i+1}\,, &\quad& \text{ for } 1\leq i\leq n-2\,;
  \end{alignedat}
\end{equation}
the last two relations are called the \emph{braid relations}. When we want to
write explicitly a permutation $\mu$ in $\sg{n}$ we will use the \emph{one
  line notation}, that is the sequence $\mu_1\mu_2\cdots\mu_n :=
\mu(1)\mu(2)\cdots\mu(n)$.

A \emph{reduced word} for an element $\mu$ of $\W$ is a decomposition
$\mu=s_1\cdots s_k$ of $\mu$ into a product of generators in $S$ of
minimal length $l(\mu)$. A \emph{(right) descent} of $\mu$ is an
element $s\in S$ such that $l(\mu s) < l(\mu)$. If $\mu$ is a
permutation, this translates into $\mu_i > \mu_{i+1}$. A \emph{recoil}
(or \emph{left descent}) of $\mu$ is a descent of $\mu^{-1}$. The sets
of left and right descents of $\mu$ are denoted respectively by
$\Rec(\mu)$ and $\Des(\mu)$. The Coxeter group $\W$ comes equipped
with three natural lattice structures. Namely $\mu<\nu$, in the
\emph{Bruhat order} (resp. \emph{(weak) left Bruhat order}, resp.
\emph{(weak) right Bruhat order}) if some reduced word for $\mu$ is a
subword (resp.  a right factor, resp. left factor) of some reduced
word for $\nu$. In type $A$, the weak Bruhat orders are the usual left
and right permutohedron.

For a subset $I$ of $S$, the parabolic subgroup $W_I$ of $W$ is the
Coxeter subgroup of $W$ generated by $I$. Left and right cosets
representatives for the quotient of $W$ by $W_I$ are given
respectively by the \emph{recoil class}
\begin{equation}
  {}^I\W=\{ \mu\in W \suchthat \Rec(\mu) \cap I = \emptyset\}
\end{equation}
and the descent class
\begin{equation}
  \W^I=\{ \mu\in W \suchthat \Des(\mu) \cap I = \emptyset\}\,.
\end{equation}

For $q$ a complex number, let $\heckeW{q}$ be the
(Iwahori-)\emph{Hecke algebra} of $W$ over the field $\CC$. This
algebra of dimension $|W|$ has a linear basis $\{T_\mu\}_{\mu\in W}$,
and its multiplication is determined by
\begin{equation}
  \begin{cases}
    T_s T_\mu = (q-1) T_{sw} + (q-1) T_w & \text{ if $s\in S$ and $\lon(s\mu) < \lon(\mu)$}\,,\\
    T_\mu T_\nu = T_{\mu\nu}         & \text{ if $\lon(\mu)+\lon(\nu) = \lon(\mu\nu)$ }\,.
  \end{cases}
\end{equation}

In particular, for any element $\mu$ of $W$, we have $T_\mu=T_{s_1}\cdots
T_{s_k}$ where $s_1,\dots,s_k$ is any reduced word for $\mu$. In fact,
$\heckeW{q}$ is the algebra generated by the elements $T_s,\, s\in S$ subject
to the same relations as the $s$ themselves, except that the quadratic
relation $s^2=1$ is replaced by:
\begin{equation}
  T_s^2=(q-1)T_s+q\,.
\end{equation}

Setting $q=1$ yields back the usual group algebra $\Wa$ of $\W$.
Similarly, the $0$-Hecke algebra $\heckeW{0}$ is obtained by setting
$q=0$ in these relations. Then, the first relation becomes
$T_s^2=-T_s$~\cite{Norton.1979,Krob_Thibon.NCSF4.1997}. In this paper,
we prefer to use another set of generators $(\pi_s)_{s\in S}$
defined by $\pi_i := T_i + 1$. They also satisfy the braid-like
relations together with the quadratic relations $\pi_s^2 = \pi_s$.
Note that the $0$-Hecke algebra is thus a monoid algebra.

\subsection{Representation theory of the $0$-Hecke algebra}
\label{section.reptheo0Hecke}

In this paper, we mostly consider \emph{right} modules over algebras.
Consequently the composition of two endomorphisms $f$ and $g$ is
denoted by $fg= g \circ f$ and their action on a vector $v$ is written
$v \act f$. Thus $g \circ f (v) = g(f(v))$ is denoted $v \act fg = (v
\act f) \act g$.

Assume now that $W$ is finite.
It is known that $\heckeW{0}$ has $2^{|S|}$ simple modules, all
one-dimensional, and naturally labeled by subsets $I$ of
$S$~\cite{Norton.1979}: following the notation
of~\cite{Krob_Thibon.NCSF4.1997}, let $\eta_I$ be the generator of the
simple $\heckeW{0}$-module $S_I$ associated with $I$ in the right
regular representation. It satisfies
\begin{equation}
\label{etaI}
\eta_I \act T_i :=
\begin{cases}
  -\eta_I & \text{if $i\in I$,}\\
  0       & \text{otherwise,}
\end{cases}
\qquad\text{or equivalently}\qquad
\eta_I \act \pi_i :=
\begin{cases}
  0       & \text{if $i\in I$,}\\
  \eta_I  & \text{otherwise.}
\end{cases}
\end{equation}

The indecomposable projective module $P_I$ associated with $S_I$ (that is such
that $S_I = P_I / \rad(P_I))$ can be described as follows: it has a basis
$\{b_\mu \suchthat \Rec(\mu) = I\}$ with the action
\begin{equation}
  b_\mu \act T_s = 
  \begin{cases}
    -b_\mu     & \text{if $s\in \Des(\mu)$,}\\
    b_{\mu s}  & \text{if $s\notin \Des(\mu)$ and $\Rec(\mu s) = I$,} \\
    0         & \text{otherwise.}
  \end{cases}
\end{equation}

\subsection{Representation theory of the $0$-Hecke algebra in type $A$}

In type $A$, it is customary to index the projective and simple
modules of $\hecke{n}{0}$ by compositions of $n$. For notational
convenience, and when there is no ambiguity, we simply identify the
subset $I$ of $S=\{1,\dots,n-1\}$ and the corresponding composition
$C(I) = (i_1,\dots,i_k)$ of $n$.

The Grothendieck rings of $\hecke{n}{0}$ are naturally isomorphic to
the dual pair of Hopf algebras of quasi-symmetric functions $\qsym$ of
Gessel~\cite{Gessel.QSym.1984} and of noncommutative symmetric
functions
$\ncsf$~\cite{Gelfand_Krob_Lascoux_Leclerc_Retakh_Thibon.1995.NCSF1}
(see~\cite{Krob_Thibon.NCSF4.1997}). The reader who is not familiar
with those should refer to these papers, as we will only recall the
required notations here.

The Hopf algebra $\qsym$ of quasi-symmetric functions has two
remarkable bases, namely the \emph{monomial basis} $(M_I)_I$ and the
\emph{fundamental basis} (also called \emph{quasi-ribbon}) $(F_I)_I$.
They are related by
\begin{equation}
  F_I = \sum_{I\fin J} M_J\,
\qquad\text{or equivalently}\qquad
  M_I = \sum_{I\fin J} (-1)^{\lon(I) - \lon(J)} F_J\,.
\end{equation}
The characteristic map $S_I \mapsto F_I$ which sends the simple
$\hecke{n}{0}$-module $S_I$ to its corresponding fundamental function
$F_I$ also sends the induction product to the product of $\qsym$ and
the restriction coproduct to the coproduct of $\qsym$.

The Hopf algebra $\ncsf$ of noncommutative symmetric
functions~\cite{Gelfand_Krob_Lascoux_Leclerc_Retakh_Thibon.1995.NCSF1}
is a noncommutative analogue of the algebra of symmetric
functions~\cite{MacDonald.SF.1995}. It has for multiplicative bases
the analogues $(\Lambda^I)_I$ of the elementary symmetric functions
$(e_\lambda)_\lambda$ and as well as the analogues $(S^I)_I$ of the
complete symmetric functions $(h_\lambda)_\lambda$. The relevant basis
in the representation theory of $\hecke{n}{0}$ is the basis of so
called \emph{ribbon Schur functions} $(R_I)_I$ which is an analogue of
skew Schur functions of ribbon shape. It is related to $(\Lambda_I)_I$
and $(S_I)_I$ by
\begin{equation}
  S_I = \sum_{I\fin J} R_J\,
  \qquad\text{and}\qquad
  \Lambda_I = \sum_{I\fin J} R_{J^c}\,.
\end{equation}
Their interpretation in representation theory goes as follows. The
complete function $S^n$ is the characteristic of the trivial module
$S_{n} \approx P_{n}$, the elementary function $\Lambda^n$ being the
characteristic of the sign module $S_{1^n} \approx P_{1^n}$. An
arbitrary indecomposable projective module $P_I$ has $R_I$ for
characteristic. Once again the map $P_I \mapsto R_I$ is an isomorphism
of Hopf algebras.

Recall that $S_I$ is the semi-simple module associated to $P_I$,
giving rise to the duality between $\calG$ and $\calK$ :
\begin{equation}
  S_I = P_I/\rad(P_I)
  \qquad\text{and}\qquad
  \langle P_I\ ,\ S_J \rangle = \delta_{I,J}\,.
\end{equation}
This translates into $\qsym$ and $\ncsf$ by setting that $(F_I)_I$ and
$(R_I)_I$ are dual bases, or equivalently that $(M_I)_I$ and $(S^I)_I$
are dual bases.

\section{The algebra $\heckeWW$}
\label{section.heckesg}

Let $(\W,S)$ be a finite Coxeter group. Its group algebra $\Wa$ and
its $0$-Hecke algebra $\heckeW{0}$ can be realized simultaneously as
operator algebras by identifying the underlying vector spaces of their
right regular representations. There are several ways to do that,
depending on which basis elements of $\heckeW{0}$ we choose to
identify with elements of $\W$. It turns out that the following
identification leads to interesting properties.

Namely, consider the plain \emph{vector space} $\kW$. On the first
hand, we identify $\kW$ with the right regular representation of the
algebra $\Wa$, i.e.: $\Wa$ acts on $\kW$ by multiplication on the
right. In type $A$, this is the usual action on positions, where an
elementary transposition $s_i$ acts on a permutation $\mu := (\mu_1,
\dots, \mu_n)$ by exchanging $\mu_i$ and $\mu_{i+1}$: $\mu\act s_i=\mu
s_i$.

On the other hand, we also identify $\kW$ with the right regular
representation of the $0$-Hecke algebra $\heckeW{0}$, i.e.:
$\heckeW{0}$ acts on the right on $\kW$ by
\begin{equation}
  \mu \act \pi_s =
  \begin{cases}
    \mu   & \text{if $\lon(\mu s) < \lon(\mu)$,}\\
    \mu s & \text{otherwise}.
  \end{cases}
\end{equation}
In type $A$, the $\pi_i:=\pi_{s_i}$'s are the \emph{elementary
  decreasing bubble sort operators}:
\begin{equation}
  \mu \act \pi_i =
  \begin{cases}
    \mu      & \text{if $\mu_i > \mu_{i+1}$,}\\
    \mu s_i & \text{otherwise}.
  \end{cases}
\end{equation}

The following easy lemma will be useful in the sequel. 
\begin{lemma}
  \label{lemma.factors}
  Let $\s,\t\in\W$. Then,
  \begin{itemize}
  \item[(a)] There exists $\t'$ such that $\s.\pi_\t=\s\t'$ with
    $\lon(\s\t')=\lon(\s)+\lon(\t')$. Furthermore, $\t'=1$ if and only if
    $\t\in\W_{\Des (\s)}$.
  \item[(b)] There exists $\s'$ such that $\s.\pi_\t=\s'\t$ with
    $\lon(\s'\t)=\lon(\s')+\lon(\t)$. Furthermore, $\s'=1$ if and only if
    $\s\in\W_{\Rec (\t)}$.
  \end{itemize}
\end{lemma}
\begin{proof}
  Applying $\pi_s$ on an element $\s$ either leaves $\s$ unchanged if
  $s\in \Des (\s)$, or extends any reduced word for $\s$ by $s$
  otherwise. (a) follows by induction; in particular $\t'$ is smaller
  than $\t$ in the Bruhat order of $\W$.
  
  (b) Since $\kW$ is the right regular representation, the linear map 
  \begin{equation}
    \Phi:
    \begin{cases}
      \kW &\rightarrow \heckeW{0}\\
      \tau &\mapsto \pi_\tau
    \end{cases}
  \end{equation}
  is a morphism of $\heckeW{0}$-module. Consequently, one has $\pi_{\s.\pi_\t}
  = \pi_\s \pi_\t$ which allows us to lift the computation to the $0$-Hecke
  monoid.  There $\s$ and $\t$ play a symmetric role, and (b) follows
  from (a) by reversion of the reduced words.
\end{proof}

\begin{definition}
  The algebra $\heckeWW$ is the subalgebra of $\End(\kW)$ generated by
  both sets of operators $\{s, \pi_s\}_{s\in S}$.
\end{definition}

By construction, the algebra $\heckeWW$ contains both $\Wa$ and
$\heckeW{0}$. In fact, it contains simultaneously all the Hecke
algebras: for any values of $q$, $\heckeW{q}$
can be realized by taking the subalgebra of $\heckeWW$ generated by
the operators:
\begin{equation}
  T_s:=(q-1) (1-\pi_s) + q s\,, \quad  \text{ for } s\in S\,.
\end{equation}
A direct calculation shows that the so-defined $T_s$ actually verifies the
Hecke relation. Reciprocally, we can recover back $\heckeWW$ by choosing for
each $s$ any two generators $T_s(q_1)$ and $T_s(q_2)$ with $q_1\ne q_2$,
because for any $q,q_1,q_2$,
\begin{equation}
  T_s(q) := \frac{q-q_1}{q_2-q_1} T_s(q_1) + \frac{q-q_2}{q_1-q_2} T_s(q_2)\,.
\end{equation}
Note that setting $T_s(1) := s$ and $T_s(0) := \pi-1$ this last equation implies
the previous one when $q_1=1$ and $q_2=0$.

Let further $\opi_s:=\pi_s s$ be the operator in $\heckeWW$ which
removes the descent $s$.  In type $A$, the $\opi_i:=\opi_{s_i}$'s are
the \emph{elementary increasing bubble sort operators}:
\begin{equation}
  \mu \act \opi_i =
  \begin{cases}
    \mu     & \text{ if } \mu_i < \mu_{i+1}\,,\\
    \mu s_i & \text{otherwise.}
  \end{cases}
\end{equation}
Since $\pi_s + \opi_s$ is a symmetrizing operator, we have the identity:
\begin{equation}
  \pi_s + \opi_s = 1 + s\,.
\end{equation}
It follows that we can alternatively take as generators for $\heckeWW$
the operators $\pi_s$'s and $\opi_s$'s.

In type $A$, the natural embedding of $\ksg{n}\otimes\ksg{m}$ in
$\ksg{n+m}$ makes $(\heckesg{n})_{n\in\NN}$ into a tower of algebras,
which contains the similar towers of algebras $(\sga{n})_{n\in\NN}$
and $(\hecke{n}{q})_{n\in\NN}$.

\subsection{Basic properties}

\begin{example}
  \label{example.HS2}
  Much of the structure of $\heckeWW$ readily appears for $\W:=\sg
  2$.  Take the natural basis $(12, 21)$ of $\ksg 2$. The matrices of
  the operators $1$, $s_1$, $\pi_1$, and $\opi_1$ are respectively:
  \begin{equation}
    \begin{pmatrix}
      1 & 0\\
      0 & 1
    \end{pmatrix}, \quad
    \begin{pmatrix}
      0 & 1\\
      1 & 0
    \end{pmatrix}, \quad
    \begin{pmatrix}
      0 & 0\\
      1 & 1
    \end{pmatrix}, \quad
    \begin{pmatrix}
      1 & 1\\
      0 & 0
    \end{pmatrix}\,.
  \end{equation}
  The algebra $\heckesg 2$ is of dimension $3$ with basis $\{1,
  s_1,\pi_1\}$ and multiplication table:
  \begin{equation}
    \begin{array}{|c||c|c|c|}
      \hline
               & 1        & s_1 & \pi_1 \\\hline\hline
      1        & 1        & s_1 & \pi_1 \\\hline
      s_1 & s_1 & 1        & \pi_1 \\\hline
      \pi_1    & \pi_1    & \opi_1 = 1 + s_1 - \pi_1 & \pi_1 \\\hline
    \end{array}
  \end{equation}
  This algebra can alternatively be described by equations. Namely, take
  $f\in\End(\k\sg 2)$ with matrix
  \begin{equation}
    \begin{pmatrix}
      f_{12,12}& f_{12,21}\\
      f_{21,12}& f_{21,21}
    \end{pmatrix}\,;
  \end{equation}
  then, the following properties are equivalent:
  \begin{itemize}
  \item $f$ belongs to $\heckesg 2$;
  \item $f_{21,21} - f_{21,12} + f_{12,21} - f_{12,12}=0$;
  \item $\begin{pmatrix}1\\-1\end{pmatrix}=12-21$ is an eigenvector of $f$;
  \item $(1-s_1)f(1+s_1) = 0$.
  \end{itemize}
  In general the relations in the parabolic subalgebra
  $\k[\pi_i,\opi_i,s_i]$ of $\heckeWW$ are:
  \begin{equation}
    \label{equation.sigmapi}
    \begin{gathered}
      \begin{aligned}
        s_i\pi_i &= \pi_i\ ,   \qquad &  s_i\opi_i &= \opi_i\ , \\
        \opi_i\pi_i   &= \pi_i\ ,   &  \pi_i\opi_i    &= \opi_i\ , \\
        \pi_is_i &= \opi_i\ ,  &  \opi_is_i &= \pi_i\ , \\
      \end{aligned}\\
      \pi_i+\opi_i=1+s_i\,.
    \end{gathered}
  \end{equation}
  In particular, any two of $\{\pi_i, \opi_i, s_i\}$ can be taken as
  generators.
\end{example}

A computer exploration suggests that in type $A$ the dimension of
$\heckesg{n}$ is given by the following sequence (sequence A000275 of the
encyclopedia of integer sequences~\cite{Sloane}):
\begin{displaymath}
  1,1,3,19,211,3651,90921,3081513,136407699,7642177651,528579161353,44237263696473,
  \dots
\end{displaymath}
These are the numbers $h_n$ of pairs $(\s,\t)$ of permutations without
common descents ($\Des(\s)\cap \Des(\t)=\emptyset$). They where first
considered by Carlitz~\cite{Carlitz.1955.BesselJ0,
  Carlitz_Scoville_Vaughan.1976.1, Carlitz_Scoville_Vaughan.1976.2} as
coefficient of the doubly exponential expansion of the inverse Bessel
function $J_0$:
\begin{equation}
  \sum_{n>=0} \frac{h_n}{n!^2}x^n = \frac{1}{J_0(\sqrt{4x})}\,.
\end{equation}
Together with
Equation~(\ref{equation.sigmapi}), this leads to state the following
\begin{theorem}
  \label{theorem.heckesg.base}
  Let $\W$ be a finite Coxeter group. A vector space basis of
  $\heckeWW$ is given by the family of operators
  \begin{equation}
    B := \left\{ \s\pi_\t \suchthat
                 (\s,\t) \in \W^2,\
                 \Des(\s) \cap \Rec(\t) = \emptyset
         \right\}\,.
  \end{equation}
  In particular the dimension of $\heckeWW$ is the number $h$ of pairs
  of elements of $W$ without common descents.
\end{theorem}


Our first approach to prove this theorem was to search for a
presentation of the algebra. In type $A$, the following relations are
easily proved to hold:
\begin{equation}
  \label{equation.rewrite}
  \begin{gathered}
    \pi_{i+1} s_{i} = \pi_{i+1}\pi_{i} +
    s_{i}s_{i+1} \pi_{i}\pi_{i+1} - \pi_{i}\pi_{i+1}\pi_{i}\ , \\
    \pi_{i} s_{i+1} = \pi_{i}\pi_{i+1} +
    s_{i+1}s_{i} \pi_{i+1}\pi_{i} - \pi_{i}\pi_{i+1}\pi_{i}\ , \\
    s_i \pi_{i+1} s_i = s_{i+1} \pi_i s_{i+1}\,,
  \end{gathered}
\end{equation}
and we conjecture that they generate all relations.
\begin{conjecture}
  A presentation of $\heckesg n$ is given by the defining relations of
  $\sg n$ and $\hecke{n}{0}$ together with the relations
  $s_i\pi_i = \pi_i$ and~(\ref{equation.rewrite}).
\end{conjecture}
Using those relations as rewriting rules yields a straightening
algorithm which rewrites any expression in the $s_i$'s and $\pi_i$'s
into a linear combination of the $\s\pi_\t$.  This algorithm seems, in
practice \emph{and} with an appropriate strategy, to always terminate.
However we have no proof of this fact; moreover this algorithm is not
efficient, due to the combinatorial explosion of the number and length
of words in intermediate results.


\begin{example}
  In the following computation for $\W=\sg 8$, we multiply some
  element of the Hecke algebra by successive elementary transposition;
  we use respectively the short hand notation $\sigma_{[154]}$ and
  $\pi_{[154]}$ for the products $s_1s_5s_4$ and $\pi_1\pi_5\pi_4$:
  \scriptsize
  \begin{equation}
    \begin{aligned}
      \pi_{[1765432]} \sigma_{[1]} &=
      \sigma_{[1234567]}\pi_{[675645342312]}
      -\sigma_{[12345]}\pi_{[7675645342312]} \\&
      +\sigma_{[234567]}\pi_{[675645342312]}
      -\sigma_{[2345]}\pi_{[7675645342312]}
      +\sigma_{[]}\pi_{[76543212]} \\
      \pi_{[1765432]} \sigma_{[13]} &=
      \sigma_{[34567234561234567]}\pi_{[67564567345623451234]}
      -\sigma_{[212345]}\pi_{[7675645342312]}\\&
      +\sigma_{[21234567]}\pi_{[675645342312]}
      -\sigma_{[3452341234567]}\pi_{[675674567345623451234]} \\&
      +\sigma_{[2345671]}\pi_{[675645342312]}
      -\sigma_{[345234567123]}\pi_{[675674567345623451234]} \\&
      +\sigma_{[3456723456712345]}\pi_{[67564567345623451234]}
      -\sigma_{[23451]}\pi_{[7675645342312]} \\&
      -\sigma_{[345672345612345]}\pi_{[767564567345623451234]}
      +\sigma_{[21]}\pi_{[765432312]} \\&
      +\sigma_{[345234123]}\pi_{[7675674567345623451234]}
      +\sigma_{[]}\pi_{[765432123]}
      -\sigma_{[]}\pi_{[7654323123]} \\
      \pi_{[1765432]} \sigma_{[135]} &=
      \sigma_{[56745634567234561234567]}\pi_{[675645673456234567123456]}+
      \text{$38$ shorter terms} \\
      \pi_{[1765432]} \sigma_{[1357]} &=
      \sigma_{[7656745634567234561234567]}\pi_{[675645673456234567123456]}
      + \text{$116$ shorter terms}
  \end{aligned}
\end{equation}
\end{example}

Encountering those difficulties does not come as a surprise. The
properties of such algebras often become clearer when considering
their concrete representations (typically as operator algebras) rather
than their abstract presentation. Here,
theorem~\ref{theorem.heckesg.base} as well as the representation
theory of $\heckeWW$ follow from the structural characterization of
$\heckeWW$ as the algebra of operators preserving certain
antisymmetries to be explained below in
Section~\ref{subsection.structuralDefinition}

\subsubsection{Variants}

As mentioned previously, the original goal of the definition of
$\heckeWW$ was to put together a Coxeter group $\W$ and its $0$-Hecke
algebra $\heckeW{0}$. Identifying their right regular representation
on the canonical basis is just one possible mean. We explore quickly
here some variants, and mention alternative constructions of
$\heckeWW$ (mostly in type $A$) in
Sections~\ref{section.alternativeDefinitions}
and~\ref{section.researchInProgress}.

A first variant is to still consider the right regular actions of $\W$
and of $\heckeW{0}$ but this time on $\W$ itself. In other words, to
consider the \emph{monoid} $\langle s, \pi_s\rangle_{s\in S}$
generated by the operators $s$ and $\pi_s$. In type $A$, the sizes of
those monoids for $n=1,2,3,4$ are $1,4,66,6264$, which are strictly
bigger than the corresponding dimensions of $\heckesg n$, in
particular because we lose the linear relations $1+\pi_s s = 1+s$.
Incidentally, an interesting question is to find a presentation of
this monoid. If instead one takes the monoid $\langle \pi_s,
\opi_s\rangle$, the sizes are $1,3,23,477$.

%
Another natural approach, in type $A$, is to start from the usual
action of $\sg n$ on the ring of polynomials $\k[x_1,\dots,x_n]$
together with the action of $\hecke n 0$ by \emph{isobaric divided
  differences} (see~\cite{Lascoux.2003.CBMS}). Note that the divided
differences being symmetrizing operators, this is in fact more a
variant on the adjoint $\heckesg n^*$ of $\heckesg n$ (see next
section). Again the obtained algebras are bigger: $1,3,20,254,\dots$,
in particular because we lose the two first relations of
Equation~(\ref{equation.rewrite}).

\subsection{$\heckeWW$ as algebra of antisymmetry-preserving operators}
\label{subsection.structuralDefinition}

Let $\ls$ be the \emph{right operator} in $\End(\kW)$ describing the
action of $s_i$ by multiplication \emph{on the left} (action on values
in type $A$). Namely $\ls$ is defined by
\begin{equation}
  \label{eq.def.lsi}
  \sigma \act \ls := s\sigma\,.
\end{equation}
A vector $v$ in $\kW$ is \emph{left $s$-symmetric} (resp.
\emph{antisymmetric}) if $v \act \ls = v$ (resp. $v \act \ls =
-v$). The subspace of left $s$-symmetric (resp. antisymmetric) vectors
can be alternatively described as the image (resp. kernel) of the
\emph{quasi-idempotent} (idempotent up to a scalar) operator
$1+\ls$, or as the kernel (resp.  image) of the quasi-idempotent
operator $1-\ls$.
\begin{theorem}
  \label{theorem.heckesg.characterization}%
  $\heckeWW$ is the subspace of $\End(\kW)$ defined by the $|S|$
  \emph{idempotent sandwich equations}:
  \begin{equation}
    (1 - \ls) \ f \ (1 + \ls) = 0\,,
    \quad  \text{ for } s\in S\,.
  \end{equation}
  In other words, $\heckeWW$ is the subalgebra of those operators in
  $\End(\kW)$ which preserve left antisymmetries.
\end{theorem}
Note that, $\ls$ being self-adjoint for the canonical scalar product
of $\kW$ (making $W$ into an orthonormal basis), the adjoint algebra
of $\heckeWW$ satisfies the equations:
\begin{equation}
  (1 + \ls) \  f\  (1 - \ls) = 0\,,
    \quad  \text{ for } s\in S\,;
\end{equation}
thus, it is the subalgebra of those operators in $\End(\kW)$ which
preserve left symmetries. Furthermore the group algebra $\Wa$ of $\W$
can be described as the subalgebra of those operators in $\End(\kW)$
which preserve both left symmetries and antisymmetries; it is
therefore the intersection of $\heckeWW$ and its adjoint $\heckeWW^*$.


\begin{proof}[Proof of theorems~\ref{theorem.heckesg.base}
  and \ref{theorem.heckesg.characterization}]%
  We proceed using three lemmas that occupy the rest of this
  subsection. We first exhibit a triangularity property of the
  operators in $B$; this proves that they are linearly independent, so
  that $\dim \heckeWW \geq h$ (lemma~\ref{lemma.independent}).  Then
  we prove that the operators in $\heckeWW$ preserve all left
  antisymmetries (lemma~\ref{lemma.antisymmetries}). Finally we
  extract from the sandwich equations $\dim \End(\kW)-h$ independent
  linear forms which are annihilated by all left antisymmetry
  preserving operators in $\End(\kW)$
  (lemma~\ref{lemma.relations.independent}).  Altogether, it follows
  simultaneously that $\heckeWW$ has dimension $h$ with $B$ as basis,
  and that $\heckeWW$ is the full subspace of left antisymmetry
  preserving operators.
\end{proof}

Let $<$ be any linear extension of the right Bruhat order on $\W$.
Given an endomorphism $f$ of $\kW$, we order the rows and columns
of its matrix $M:=[f_{\mu\nu}]$ according to $<$, and denote by
$\init(f) := \min\{ \mu \suchthat \exists \nu, f_{\mu\nu}\ne 0\}$ the index of
the first non zero row of $M$.
\begin{lemma}
  \label{lemma.independent}
  (a) Let $f:=\s \pi_\t$ in $B$. Then, $\init(f)=\t$, and
  \begin{equation}
    f_{\t\nu} =
    \begin{cases}
      1 & \text{ if } \nu\in \W_{\Rec(\t)} \s^{-1}\,,\\
      0 & \text{ otherwise.}
    \end{cases}
  \end{equation}

  (b) The family $B$ is free.
\end{lemma}

\begin{proof}
  (a) is a direct corollary of Lemma~\ref{lemma.factors}~(b).

  (b) follows by triangularity: the operator $\s \pi_\t$ has
  coefficient $m_{\t\s^{-1}}=1$, whereas any other operator
  $\s'\pi_\t'$ such that $D(\s')\cap R(\t) = \emptyset$, $\t'\leq \t$
  and $\s'\ne\s$ has coefficient $m_{\t\s^{-1}}=0$.
\end{proof}

\begin{lemma}
  \label{lemma.antisymmetries}
  The operators in $\heckeWW$ preserve all left antisymmetries.
\end{lemma}
\begin{proof}
  It is sufficient to prove that the generators $s$ and $\pi_s$ of
  $\heckeWW$ preserve any left antisymmetry. For a generator $s$, this
  is obvious since the actions of $\ls$ and $s$ commute.  Let now $v$
  be an $s'$-antisymmetric vector; without loss of generality, we may
  assume that $v=(1-s')\sigma$ where $\sigma$ is some permutation
  without recoil at position $s'$. We use the same linear isomorphism
  $\Phi$ as in lemma~\ref{lemma.factors} to lift the computation to
  the $0$-Hecke algebra and use its associativity:
  \begin{equation}
    \begin{aligned}
      v\act\pi_s 
      & =\Phi^{-1}(((1-\pi_{s'})\pi_\sigma)\pi_s) \\
      & = \Phi^{-1}((1-\pi_{s'}) (\pi_\s\pi_s))       \\
      & = \Phi^{-1}((1-\pi_{s'}) \pi_{\s\ \act \pi_s}) \\
      & =
      \begin{cases}
        0         & \text{ if } s\in \Rec(\s \act \pi_s) \,,\\
        (1-s') (\s \act \pi_s) & \text{ otherwise. }
      \end{cases}
    \end{aligned}
  \end{equation}
  Therefore, $v \act \pi_s$ is again $s'$-antisymmetric.
\end{proof}

We now turn to the explicit description of the sandwich equations.
Given an endomorphism $f$ of $\kW$, denote by
$(f_{\mu,\nu})_{\mu,\nu}$ the coefficients of its matrix in the
natural group basis. Given two elements $\mu,\nu$ in $\W$ and a simple
reflection $s\in S$, let $R_{\mu,\nu,i}$ be the linear form on
$\End(\kW)$ which computes the $(\mu,\nu)$ coefficient of the matrix
of $(1-\ls)f(1+\ls)$:
\begin{equation}
  R_{\mu,\nu,s}(f) := ((1-\ls)f(1+\ls))_{\mu,\nu}
  = f_{\mu,\nu} - f_{\mu,s\nu} + f_{s\mu,\nu} - f_{s\mu,s_i\nu}
  \,.
\end{equation}
By construction, $R_{\mu,\nu,i}$ annihilates any operator which
preserves $s$-antisymmetries.

Given a pair $(\mu,\nu)$ of elements of $\W$ having at least one
descent in common, set $R_{\mu,\nu}=R_{\mu,\nu,s}$, where $s$ is the
smallest common descent of $\mu$ and $\nu$ (the choice of the common
descent $s$ is, in fact, irrelevant for our purposes). Finally, let
\begin{equation}
  R:=\{ R_{\mu,\nu} \suchthat \Des(\mu) \cap \Des(\nu) \ne \emptyset \}\,.
\end{equation}
For example in type $A$ we have:
\begin{equation}
  \begin{aligned}
    \text{for $n=1$: } & R=\{\}\\
    \text{for $n=2$: } & R=\{ R_{21,21,1} \} = \{ f_{21,21} - f_{21,12} + f_{12,21} - f_{12,12}\}\\
    \text{for $n=3$: } & R=\{ \underbrace{R_{213,213,1}, R_{213,312,1}, R_{213,321},
    \dots, R_{321,321,1}}_{\text{$17=6^2-19$ items}} \} \\
  \end{aligned}
\end{equation}
For $n=2$, the reader will recognize the linear relation described in
example~\ref{example.HS2}.

\begin{lemma}
  \label{lemma.relations.independent}
  The $|\W|^2-h$ linear forms in $R$ are linearly independent.
\end{lemma}
\begin{proof}
  Take some linear form $R_{\mu,\nu}=R_{\mu,\nu,i}$ in $R$, and
  represent it as the $|\W|\times|\W|$ array of its values on the
  elements of the canonical basis of $\End(\W)$, with the rows and
  columns sorted as previously. For example, here is the array for
  $R_{213, 312}=R_{213,312,1}$ in type $A_2$:
  \begin{displaymath}
    \begin{array}{cccccc|c}
      123 & 132 & 213 & 231 & 312 & 321 & \\\hline
      0  & -1  &  0  &  0  &  1  &  0  & 123 \\
      0  &  0  &  0  &  0  &  0  &  0  & 132 \\
      0  & -1  &  0  &  0  &  1  &  0  & 213 \\
      0  &  0  &  0  &  0  &  0  &  0  & 231 \\
      0  &  0  &  0  &  0  &  0  &  0  & 312 \\
      0  &  0  &  0  &  0  &  0  &  0  & 321 \\
    \end{array}
  \end{displaymath}
  This array has a coefficient $1$ at position $(\mu,\nu)$.  Since $s$
  is both a descent of $\mu$ and $\nu$, $s\mu < \mu$, and $s\nu <
  \nu$; so the three other non zero coefficients are either strictly
  higher or strictly to the left in the array.  Furthermore, no other
  linear form in $R$ has a non-zero coefficient at position
  $(\mu,\nu)$.  Hence, by triangularity the linear forms in $R$ are
  linearly independent.
\end{proof}

\subsection{Representation theory}

Due to the particular structure of $\heckeWW$ as an operator algebra,
the easiest thing to start with is the study of projective modules.
Along the way we define a particular basis of $\kW$ which plays a key
role for the representation theory.

\subsubsection{Projective modules}

Recall that $\heckeWW$ is the algebra of operators preserving left
antisymmetries. Thus, given $I\subset S$, it is natural to introduce
the $\heckeWW$-submodule
\begin{equation}
  \label{eq.def.P.I}
  P_I := \bigcap_{s \in I} \ker (1+\ls)\, .
\end{equation}
of the vectors in $\kW$ which are $s$-antisymmetric for all $s\in I$.
For example, $P_S$ is one dimensional, and spanned by $\sum_{\nu\in \W}
(-1)^{\lon(\nu)} \nu$, whereas $P_\emptyset=\kW$.

The goal of this section is to prove that the family of modules
$(P_I)_{I\subset S}$ forms a complete set of representatives of the
indecomposable projective modules of $\heckeWW$. First, we need a more
practical definition of $P_I$.
\begin{lemma}
  Let $I\subset S$. Then $P_I$ is the $\heckeWW$-submodule of $\kW$
  generated by
  \begin{equation}
    v_I := \sum_{\nu \in \W_I} (-1)^{\lon(\nu)} \nu\,,
  \end{equation}
  or equivalently the $\Wa$-submodule generated by $v_I$.
\end{lemma}
\begin{proof}
  First, it is clear that $v_I$ belongs to $P_I$. Since $\kW$ is the
  right regular representation of $\W$, we may temporarily identify
  $\kW$ and $\Wa$. There, it is well known that $v_I$ is an idempotent
  (up to a scalar factor). Take in general $u\in P_I$. For any $s$ in
  $I$, one has $(1+s)u=0$, that is $su=-u$. It follows that $v_I u =
  |W_I|u$, and we can conclude that $P_I = v_I\act\Wa$.
\end{proof}
Actually, as we will see later, it is also an idempotent in the $0$-Hecke
algebra and even in the generic Hecke algebra.

\begin{figure}
  \subfigure[Projective modules $P_{\{\}}$ and $P_{\{1\}}$ of $\heckesg 2$]{
    \centerline{
      \begin{tikzpicture}[scale=.5,>=latex,join=bevel,]
\pgfsetlinewidth{.2pt}
\tiny%
  \definecolor{newcol}{rgb}{0.0,0.0,0.0};
  \pgfsetcolor{newcol}
  \draw [->] (33bp,94bp) .. controls (28bp,89bp) and (22bp,83bp)  .. (19bp,76bp) .. controls (15bp,63bp) and (15bp,58bp)  .. (19bp,46bp) .. controls (21bp,42bp) and (23bp,37bp)  .. (33bp,26bp);
  \draw (35bp,61bp) node {$\begin{array}{r}s_1\\\pi_1\\\end{array}$};
  \draw [->] (43bp,26bp) .. controls (44bp,31bp) and (44bp,39bp)  .. (45bp,46bp) .. controls (45bp,59bp) and (45bp,62bp)  .. (45bp,76bp) .. controls (45bp,79bp) and (44bp,81bp)  .. (43bp,94bp);
  \draw (57bp,61bp) node {$\begin{array}{r}s_1\\\end{array}$};
  \draw [->] (46bp,26bp) .. controls (57bp,40bp) and (73bp,37bp)  .. (73bp,19bp) .. controls (73bp,5bp) and (63bp,0bp)  .. (46bp,12bp);
  \draw (86bp,19bp) node {$\begin{array}{r}\pi_1\\\end{array}$};
\begin{scope}
  \pgfsetstrokecolor{black}
  \draw (42bp,101bp) node {$12$};
\end{scope}
\begin{scope}
  \pgfsetstrokecolor{black}
  \draw (42bp,19bp) node {$21$};
\end{scope}
\end{tikzpicture} \qquad
      \begin{tikzpicture}[scale=.5,>=latex,join=bevel,]
\pgfsetlinewidth{.2pt}
\tiny%
  \definecolor{newcol}{rgb}{0.0,0.0,0.0};
  \pgfsetcolor{newcol}
  \draw [->] (26bp,14bp) .. controls (35bp,15bp) and (44bp,14bp)  .. (44bp,9bp) .. controls (44bp,6bp) and (41bp,4bp)  .. (26bp,4bp);
  \draw (60bp,9bp) node {$\begin{array}{r}- s_1\\\end{array}$};
\begin{scope}
  \pgfsetstrokecolor{black}
  \draw (13bp,9bp) node {$12$};
\end{scope}
\end{tikzpicture}}}

  \subfigure[Projective modules $P_{\{\}}$, $P_{\{1\}}$, $P_{\{2\}}$, and  $P_{\{1,2\}}$ of $\heckesg 3$]{
    \centerline{
      \begin{tikzpicture}[scale=.5,>=latex,join=bevel,]
\pgfsetlinewidth{.2pt}
\tiny%
\begin{scope}
  \pgfsetstrokecolor{black}
  \pgfsetfillcolor{white}
  \pgfsetstrokecolor{white}
  \filldraw (-2bp,-2bp) -- (-2bp,2bp) -- (2bp,2bp) -- (2bp,-2bp) -- cycle;
\end{scope}
\begin{scope}
  \pgfsetstrokecolor{black}
  \pgfsetfillcolor{white}
  \pgfsetstrokecolor{white}
  \filldraw (-2bp,-2bp) -- (-2bp,2bp) -- (2bp,2bp) -- (2bp,-2bp) -- cycle;
\end{scope}
  \definecolor{newcol}{rgb}{0.0,0.0,0.0};
  \pgfsetcolor{newcol}
  \draw [->] (81bp,275bp) .. controls (62bp,270bp) and (32bp,262bp)  .. (25bp,253bp) .. controls (17bp,241bp) and (22bp,224bp)  .. (31bp,203bp);
  \draw (41bp,238bp) node {$\begin{array}{r}s_1\\\pi_1\\\end{array}$};
  \draw [->] (98bp,271bp) .. controls (102bp,260bp) and (110bp,239bp)  .. (121bp,223bp) .. controls (123bp,219bp) and (126bp,215bp)  .. (136bp,203bp);
  \draw (134bp,238bp) node {$\begin{array}{r}s_2\\\pi_2\\\end{array}$};
  \draw [->] (41bp,203bp) .. controls (47bp,214bp) and (61bp,235bp)  .. (75bp,253bp) .. controls (77bp,256bp) and (80bp,260bp)  .. (90bp,271bp);
  \draw (87bp,238bp) node {$\begin{array}{r}s_1\\\end{array}$};
  \draw [->] (41bp,203bp) .. controls (52bp,217bp) and (69bp,214bp)  .. (69bp,196bp) .. controls (69bp,182bp) and (58bp,177bp)  .. (41bp,189bp);
  \draw (82bp,196bp) node {$\begin{array}{r}\pi_1\\\end{array}$};
  \draw [->] (30bp,189bp) .. controls (26bp,184bp) and (21bp,176bp)  .. (19bp,169bp) .. controls (16bp,156bp) and (14bp,151bp)  .. (19bp,139bp) .. controls (21bp,134bp) and (26bp,129bp)  .. (38bp,119bp);
  \draw (35bp,154bp) node {$\begin{array}{r}s_2\\\pi_2\\\end{array}$};
  \draw [->] (139bp,189bp) .. controls (136bp,184bp) and (132bp,176bp)  .. (130bp,169bp) .. controls (126bp,156bp) and (125bp,151bp)  .. (130bp,139bp) .. controls (131bp,135bp) and (133bp,131bp)  .. (142bp,119bp);
  \draw (143bp,154bp) node {$\begin{array}{r}s_1\\\pi_1\\\end{array}$};
  \draw [->] (147bp,203bp) .. controls (151bp,215bp) and (156bp,237bp)  .. (147bp,253bp) .. controls (141bp,262bp) and (131bp,268bp)  .. (111bp,275bp);
  \draw (163bp,238bp) node {$\begin{array}{r}s_2\\\end{array}$};
  \draw [->] (149bp,203bp) .. controls (160bp,217bp) and (177bp,214bp)  .. (177bp,196bp) .. controls (177bp,182bp) and (166bp,177bp)  .. (149bp,189bp);
  \draw (190bp,196bp) node {$\begin{array}{r}\pi_2\\\end{array}$};
  \draw [->] (48bp,105bp) .. controls (44bp,93bp) and (38bp,70bp)  .. (47bp,55bp) .. controls (57bp,39bp) and (77bp,30bp)  .. (103bp,24bp);
  \draw (60bp,70bp) node {$\begin{array}{r}s_1\\\pi_1\\\end{array}$};
  \draw [->] (51bp,119bp) .. controls (50bp,130bp) and (49bp,151bp)  .. (45bp,169bp) .. controls (44bp,173bp) and (43bp,176bp)  .. (39bp,189bp);
  \draw (60bp,154bp) node {$\begin{array}{r}s_2\\\end{array}$};
  \draw [->] (56bp,119bp) .. controls (67bp,133bp) and (84bp,130bp)  .. (84bp,112bp) .. controls (84bp,98bp) and (73bp,93bp)  .. (56bp,105bp);
  \draw (97bp,112bp) node {$\begin{array}{r}\pi_2\\\end{array}$};
  \draw [->] (151bp,119bp) .. controls (153bp,124bp) and (155bp,132bp)  .. (156bp,139bp) .. controls (157bp,152bp) and (158bp,155bp)  .. (156bp,169bp) .. controls (155bp,173bp) and (154bp,176bp)  .. (148bp,189bp);
  \draw (168bp,154bp) node {$\begin{array}{r}s_1\\\end{array}$};
  \draw [->] (154bp,119bp) .. controls (165bp,133bp) and (182bp,130bp)  .. (182bp,112bp) .. controls (182bp,98bp) and (171bp,93bp)  .. (154bp,105bp);
  \draw (195bp,112bp) node {$\begin{array}{r}\pi_1\\\end{array}$};
  \draw [->] (147bp,105bp) .. controls (142bp,91bp) and (130bp,59bp)  .. (120bp,29bp);
  \draw (151bp,70bp) node {$\begin{array}{r}s_2\\\pi_2\\\end{array}$};
  \draw [->] (113bp,29bp) .. controls (104bp,41bp) and (87bp,65bp)  .. (73bp,85bp) .. controls (70bp,89bp) and (66bp,93bp)  .. (57bp,105bp);
  \draw (104bp,70bp) node {$\begin{array}{r}s_1\\\end{array}$};
  \draw [->] (122bp,29bp) .. controls (134bp,46bp) and (151bp,43bp)  .. (151bp,22bp) .. controls (151bp,5bp) and (140bp,0bp)  .. (122bp,15bp);
  \draw (164bp,22bp) node {$\begin{array}{r}\pi_1\\\pi_2\\\end{array}$};
  \draw [->] (133bp,27bp) .. controls (145bp,32bp) and (160bp,41bp)  .. (168bp,55bp) .. controls (174bp,66bp) and (172bp,72bp)  .. (168bp,85bp) .. controls (166bp,89bp) and (164bp,93bp)  .. (155bp,105bp);
  \draw (183bp,70bp) node {$\begin{array}{r}s_2\\\end{array}$};
\begin{scope}
  \pgfsetstrokecolor{black}
  \draw (96bp,278bp) node {$123$};
\end{scope}
\begin{scope}
  \pgfsetstrokecolor{black}
  \draw (118bp,22bp) node {$321$};
\end{scope}
\begin{scope}
  \pgfsetstrokecolor{black}
  \draw (36bp,196bp) node {$213$};
\end{scope}
\begin{scope}
  \pgfsetstrokecolor{black}
  \draw (144bp,196bp) node {$132$};
\end{scope}
\begin{scope}
  \pgfsetstrokecolor{black}
  \draw (51bp,112bp) node {$231$};
\end{scope}
\begin{scope}
  \pgfsetstrokecolor{black}
  \draw (149bp,112bp) node {$312$};
\end{scope}
\end{tikzpicture} \qquad
      \begin{tikzpicture}[scale=.5,>=latex,join=bevel,]
\pgfsetlinewidth{.2pt}
\tiny%
  \definecolor{newcol}{rgb}{0.0,0.0,0.0};
  \pgfsetcolor{newcol}
  \draw [->] (47bp,203bp) .. controls (58bp,217bp) and (75bp,214bp)  .. (75bp,196bp) .. controls (75bp,182bp) and (64bp,177bp)  .. (47bp,189bp);
  \draw (91bp,196bp) node {$\begin{array}{r}- s_1\\\end{array}$};
  \draw [->] (34bp,189bp) .. controls (28bp,184bp) and (22bp,176bp)  .. (19bp,169bp) .. controls (15bp,156bp) and (15bp,151bp)  .. (19bp,139bp) .. controls (21bp,135bp) and (23bp,130bp)  .. (33bp,119bp);
  \draw (35bp,154bp) node {$\begin{array}{r}s_2\\\pi_2\\\end{array}$};
  \draw [->] (34bp,105bp) .. controls (28bp,100bp) and (22bp,92bp)  .. (19bp,85bp) .. controls (15bp,72bp) and (16bp,67bp)  .. (19bp,55bp) .. controls (21bp,49bp) and (25bp,43bp)  .. (35bp,29bp);
  \draw (35bp,70bp) node {$\begin{array}{r}s_1\\\pi_1\\\end{array}$};
  \draw [->] (43bp,119bp) .. controls (44bp,124bp) and (44bp,132bp)  .. (45bp,139bp) .. controls (45bp,152bp) and (45bp,155bp)  .. (45bp,169bp) .. controls (45bp,172bp) and (44bp,176bp)  .. (43bp,189bp);
  \draw (57bp,154bp) node {$\begin{array}{r}s_2\\\end{array}$};
  \draw [->] (47bp,119bp) .. controls (58bp,133bp) and (75bp,130bp)  .. (75bp,112bp) .. controls (75bp,98bp) and (64bp,93bp)  .. (47bp,105bp);
  \draw (88bp,112bp) node {$\begin{array}{r}\pi_2\\\end{array}$};
  \draw [->] (43bp,29bp) .. controls (43bp,36bp) and (44bp,46bp)  .. (45bp,55bp) .. controls (45bp,68bp) and (45bp,71bp)  .. (45bp,85bp) .. controls (45bp,88bp) and (44bp,92bp)  .. (43bp,105bp);
  \draw (57bp,70bp) node {$\begin{array}{r}s_1\\\end{array}$};
  \draw [->] (46bp,29bp) .. controls (58bp,46bp) and (75bp,43bp)  .. (75bp,22bp) .. controls (75bp,5bp) and (64bp,0bp)  .. (46bp,15bp);
  \draw (91bp,22bp) node {$\begin{array}{r}\pi_1\\- s_2\\\end{array}$};
\begin{scope}
  \pgfsetstrokecolor{black}
  \draw (42bp,196bp) node {$123$};
\end{scope}
\begin{scope}
  \pgfsetstrokecolor{black}
  \draw (42bp,112bp) node {$132$};
\end{scope}
\begin{scope}
  \pgfsetstrokecolor{black}
  \draw (42bp,22bp) node {$312$};
\end{scope}
\end{tikzpicture} \quad
      \begin{tikzpicture}[scale=.5,>=latex,join=bevel,]
\pgfsetlinewidth{.2pt}
\tiny%
  \definecolor{newcol}{rgb}{0.0,0.0,0.0};
  \pgfsetcolor{newcol}
  \draw [->] (34bp,189bp) .. controls (28bp,184bp) and (22bp,176bp)  .. (19bp,169bp) .. controls (15bp,156bp) and (15bp,151bp)  .. (19bp,139bp) .. controls (21bp,135bp) and (23bp,130bp)  .. (33bp,119bp);
  \draw (35bp,154bp) node {$\begin{array}{r}s_1\\\pi_1\\\end{array}$};
  \draw [->] (47bp,203bp) .. controls (58bp,217bp) and (75bp,214bp)  .. (75bp,196bp) .. controls (75bp,182bp) and (64bp,177bp)  .. (47bp,189bp);
  \draw (91bp,196bp) node {$\begin{array}{r}- s_2\\\end{array}$};
  \draw [->] (43bp,119bp) .. controls (44bp,124bp) and (44bp,132bp)  .. (45bp,139bp) .. controls (45bp,152bp) and (45bp,155bp)  .. (45bp,169bp) .. controls (45bp,172bp) and (44bp,176bp)  .. (43bp,189bp);
  \draw (57bp,154bp) node {$\begin{array}{r}s_1\\\end{array}$};
  \draw [->] (47bp,119bp) .. controls (58bp,133bp) and (75bp,130bp)  .. (75bp,112bp) .. controls (75bp,98bp) and (64bp,93bp)  .. (47bp,105bp);
  \draw (88bp,112bp) node {$\begin{array}{r}\pi_1\\\end{array}$};
  \draw [->] (34bp,105bp) .. controls (28bp,100bp) and (22bp,92bp)  .. (19bp,85bp) .. controls (15bp,72bp) and (16bp,67bp)  .. (19bp,55bp) .. controls (21bp,49bp) and (25bp,43bp)  .. (35bp,29bp);
  \draw (35bp,70bp) node {$\begin{array}{r}s_2\\\pi_2\\\end{array}$};
  \draw [->] (46bp,29bp) .. controls (58bp,46bp) and (75bp,43bp)  .. (75bp,22bp) .. controls (75bp,5bp) and (64bp,0bp)  .. (46bp,15bp);
  \draw (91bp,22bp) node {$\begin{array}{r}- s_1\\\pi_2\\\end{array}$};
  \draw [->] (43bp,29bp) .. controls (43bp,36bp) and (44bp,46bp)  .. (45bp,55bp) .. controls (45bp,68bp) and (45bp,71bp)  .. (45bp,85bp) .. controls (45bp,88bp) and (44bp,92bp)  .. (43bp,105bp);
  \draw (57bp,70bp) node {$\begin{array}{r}s_2\\\end{array}$};
\begin{scope}
  \pgfsetstrokecolor{black}
  \draw (42bp,196bp) node {$123$};
\end{scope}
\begin{scope}
  \pgfsetstrokecolor{black}
  \draw (42bp,112bp) node {$213$};
\end{scope}
\begin{scope}
  \pgfsetstrokecolor{black}
  \draw (42bp,22bp) node {$231$};
\end{scope}
\end{tikzpicture} \quad
      \begin{tikzpicture}[scale=.5,>=latex,join=bevel,]
\pgfsetlinewidth{.2pt}
\tiny%
  \definecolor{newcol}{rgb}{0.0,0.0,0.0};
  \pgfsetcolor{newcol}
  \draw [->] (30bp,20bp) .. controls (39bp,21bp) and (48bp,19bp)  .. (48bp,15bp) .. controls (48bp,13bp) and (45bp,11bp)  .. (30bp,10bp);
  \draw (64bp,15bp) node {$\begin{array}{r}- s_1\\- s_2\\\end{array}$};
\begin{scope}
  \pgfsetstrokecolor{black}
  \draw (15bp,15bp) node {$123$};
\end{scope}
\end{tikzpicture}}}

  \subfigure[Projective module $P_{\{1,3\}}$ of $\heckesg 4$]{
    \centerline{
      \begin{tikzpicture}[scale=.5,>=latex,join=bevel,]
\pgfsetlinewidth{.2pt}
\tiny%
\begin{scope}
  \pgfsetstrokecolor{black}
  \pgfsetfillcolor{white}
  \pgfsetstrokecolor{white}
  \filldraw (-2bp,-2bp) -- (-2bp,2bp) -- (2bp,2bp) -- (2bp,-2bp) -- cycle;
\end{scope}
  \definecolor{newcol}{rgb}{0.0,0.0,0.0};
  \pgfsetcolor{newcol}
  \draw [->] (94bp,411bp) .. controls (107bp,428bp) and (126bp,425bp)  .. (126bp,404bp) .. controls (126bp,387bp) and (113bp,382bp)  .. (94bp,397bp);
  \draw (142bp,404bp) node {$\begin{array}{r}- s_1\\- s_3\\\end{array}$};
  \draw [->] (83bp,397bp) .. controls (78bp,391bp) and (71bp,381bp)  .. (67bp,371bp) .. controls (64bp,358bp) and (63bp,353bp)  .. (67bp,341bp) .. controls (69bp,337bp) and (71bp,332bp)  .. (81bp,321bp);
  \draw (83bp,356bp) node {$\begin{array}{r}s_2\\\pi_2\\\end{array}$};
  \draw [->] (72bp,309bp) .. controls (53bp,304bp) and (25bp,295bp)  .. (19bp,287bp) .. controls (9bp,273bp) and (18bp,253bp)  .. (32bp,231bp);
  \draw (35bp,272bp) node {$\begin{array}{r}s_1\\\pi_1\\\end{array}$};
  \draw [->] (91bp,321bp) .. controls (92bp,326bp) and (92bp,334bp)  .. (93bp,341bp) .. controls (93bp,354bp) and (93bp,357bp)  .. (93bp,371bp) .. controls (92bp,376bp) and (92bp,382bp)  .. (91bp,397bp);
  \draw (105bp,356bp) node {$\begin{array}{r}s_2\\\end{array}$};
  \draw [->] (95bp,321bp) .. controls (108bp,335bp) and (126bp,332bp)  .. (126bp,314bp) .. controls (126bp,300bp) and (114bp,295bp)  .. (95bp,307bp);
  \draw (139bp,314bp) node {$\begin{array}{r}\pi_2\\\end{array}$};
  \draw [->] (92bp,307bp) .. controls (95bp,295bp) and (103bp,273bp)  .. (115bp,257bp) .. controls (120bp,250bp) and (127bp,243bp)  .. (142bp,231bp);
  \draw (128bp,272bp) node {$\begin{array}{r}s_3\\\pi_3\\\end{array}$};
  \draw [->] (41bp,231bp) .. controls (46bp,243bp) and (56bp,267bp)  .. (68bp,287bp) .. controls (71bp,291bp) and (74bp,295bp)  .. (84bp,307bp);
  \draw (83bp,272bp) node {$\begin{array}{r}s_1\\\end{array}$};
  \draw [->] (42bp,231bp) .. controls (55bp,248bp) and (74bp,245bp)  .. (74bp,224bp) .. controls (74bp,207bp) and (61bp,202bp)  .. (42bp,217bp);
  \draw (90bp,224bp) node {$\begin{array}{r}\pi_1\\- s_2\\\end{array}$};
  \draw [->] (33bp,217bp) .. controls (26bp,204bp) and (13bp,179bp)  .. (24bp,161bp) .. controls (34bp,146bp) and (51bp,137bp)  .. (77bp,131bp);
  \draw (37bp,176bp) node {$\begin{array}{r}s_3\\\pi_3\\\end{array}$};
  \draw [->] (143bp,217bp) .. controls (136bp,211bp) and (126bp,201bp)  .. (120bp,191bp) .. controls (111bp,176bp) and (104bp,158bp)  .. (97bp,135bp);
  \draw (133bp,176bp) node {$\begin{array}{r}s_1\\\pi_1\\\end{array}$};
  \draw [->] (156bp,231bp) .. controls (169bp,248bp) and (188bp,245bp)  .. (188bp,224bp) .. controls (188bp,207bp) and (175bp,202bp)  .. (156bp,217bp);
  \draw (204bp,224bp) node {$\begin{array}{r}- s_2\\\pi_3\\\end{array}$};
  \draw [->] (153bp,231bp) .. controls (153bp,244bp) and (153bp,270bp)  .. (141bp,287bp) .. controls (135bp,294bp) and (126bp,300bp)  .. (108bp,309bp);
  \draw (162bp,272bp) node {$\begin{array}{r}s_3\\\end{array}$};
  \draw [->] (113bp,135bp) .. controls (124bp,140bp) and (138bp,148bp)  .. (146bp,161bp) .. controls (154bp,175bp) and (155bp,193bp)  .. (153bp,217bp);
  \draw (165bp,176bp) node {$\begin{array}{r}s_1\\\end{array}$};
  \draw [->] (99bp,135bp) .. controls (112bp,152bp) and (131bp,149bp)  .. (131bp,128bp) .. controls (131bp,111bp) and (118bp,106bp)  .. (99bp,121bp);
  \draw (144bp,128bp) node {$\begin{array}{r}\pi_1\\\pi_3\\\end{array}$};
  \draw [->] (88bp,121bp) .. controls (83bp,115bp) and (76bp,105bp)  .. (72bp,95bp) .. controls (69bp,82bp) and (69bp,77bp)  .. (72bp,65bp) .. controls (75bp,57bp) and (79bp,48bp)  .. (89bp,33bp);
  \draw (88bp,80bp) node {$\begin{array}{r}s_2\\\pi_2\\\end{array}$};
  \draw [->] (91bp,135bp) .. controls (82bp,151bp) and (60bp,187bp)  .. (42bp,217bp);
  \draw (88bp,176bp) node {$\begin{array}{r}s_3\\\end{array}$};
  \draw [->] (99bp,33bp) .. controls (112bp,53bp) and (131bp,50bp)  .. (131bp,26bp) .. controls (131bp,6bp) and (117bp,0bp)  .. (99bp,19bp);
  \draw (147bp,26bp) node {$\begin{array}{r}- s_1\\\pi_2\\- s_3\\\end{array}$};
  \draw [->] (96bp,33bp) .. controls (96bp,41bp) and (97bp,53bp)  .. (98bp,65bp) .. controls (98bp,78bp) and (98bp,81bp)  .. (98bp,95bp) .. controls (97bp,100bp) and (97bp,106bp)  .. (96bp,121bp);
  \draw (110bp,80bp) node {$\begin{array}{r}s_2\\\end{array}$};
\begin{scope}
  \pgfsetstrokecolor{black}
  \draw (90bp,404bp) node {$1234$};
\end{scope}
\begin{scope}
  \pgfsetstrokecolor{black}
  \draw (90bp,314bp) node {$1324$};
\end{scope}
\begin{scope}
  \pgfsetstrokecolor{black}
  \draw (95bp,128bp) node {$3142$};
\end{scope}
\begin{scope}
  \pgfsetstrokecolor{black}
  \draw (95bp,26bp) node {$3412$};
\end{scope}
\begin{scope}
  \pgfsetstrokecolor{black}
  \draw (38bp,224bp) node {$3124$};
\end{scope}
\begin{scope}
  \pgfsetstrokecolor{black}
  \draw (152bp,224bp) node {$1342$};
\end{scope}
\end{tikzpicture}
    }}

  \caption{Some projective modules of $\heckesg n$, described on their
    $\{v_I \act \s \suchthat \s \in {}^I\W\}$ bases. Note that those
    are not quite combinatorial modules, as coefficients $-1$ or $0$
    may occur; in the later case, the corresponding edges are not
    drawn.  Those pictures have been produced automatically, using
    $\texttt{MuPAD-Combinat}$, \texttt{graphviz}, and
    \texttt{dot2tex}.}
  \label{figure.almostCombinatorialModule}
\end{figure}
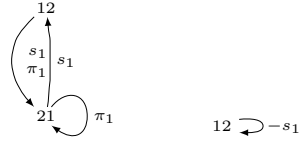
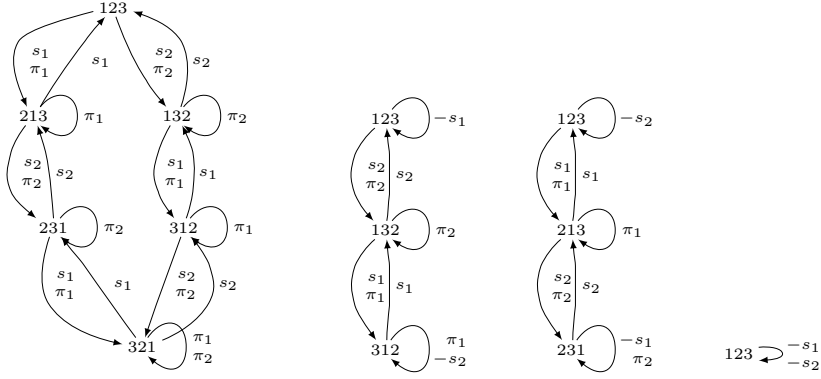
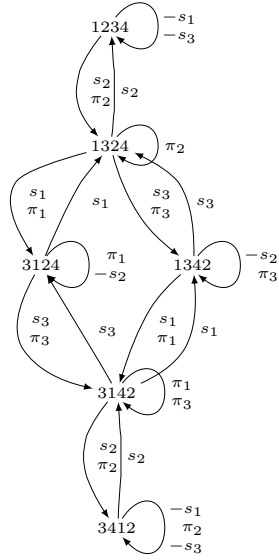

For each $\s\in W$, define $v_\s:=v_{S\backslash\Rec(\s)}\act \s$.
Note that $\s$ is the element of minimal length appearing in $v_\s$.
By triangularity, it follows that the family $(v_\s)_{\s\in\sg{n}}$
forms a vector space basis of $\kW$. See
Figures~\ref{figure.almostCombinatorialModule} and~\ref{figure.vsigma}
for examples.

The usefulness of this basis comes from the fact that it is compatible
with the module structure.
\begin{proposition}
  \label{proposition.base.PI}
  For any $I\subset S$, the module $P_I$ is of dimension $|{}^I\W|$,
  and
  \begin{equation}
    \{v_I \act \s \suchthat \s \in {}^I\W\} 
    \qquad\text{and}\qquad
    \{v_\s\       \suchthat \s \in {}^I\W\}
  \end{equation}
  are both vector space bases of $P_I$.
\end{proposition}

\begin{proof}
  First note that, by the same triangularity argument, the first
  family is free as well. Furthermore, by the previous lemma, $P_I$ is
  spanned by all the vectors $v_I\act \sigma$ with $\sigma \in W$.
  Take $\sigma\notin {}^I\W$; then $\sigma$ is of the form $\sigma =
  s\sigma'$ with $s \in I$ and $\lon(\sigma')<\lon(\sigma)$, and it
  follows that $v_I\act\sigma = -v_I\act\sigma'$. Applying induction
  on the length yields that the first family is a basis. By dimension
  count, the second family (each element $v_\s$ of which is in
  $P_{S\backslash\Rec{\s}}\subset P_I$) is also a basis of $P_I$.
\end{proof}

\begin{corollaries}
  $P_I$ is generated by $v_I$, either as $\heckeWW$, $\W$, or
  $\heckeW{0}$-module.

  $P_J\subset P_I$ if and only if $I\subset J$.
\end{corollaries}
\begin{proof}
  If a finite Coxeter group, the sets $\{w\ \suchthat\ \Des(w) = I\}$ 
  are never empty;
  therefore $P_J\subsetneq P_I$ whenever $I\subsetneq J$.
\end{proof}

In type $A$, the recoil class ${}^I\W$ is the set of the shuffles of
the words $1\cdots a_1$, $a_1+1\cdots a_2$, \dots, where $I$ is the
set $\{a_1 < a_2 < \cdots\}$. As a consequence
\begin{corollary}
  In type $A$ the module $P_I$ is of dimension $\frac{n!}{i_1!i_2!\dots
    i_k!}$, where $(i_1,\dots,i_k)$ is the composition associated to $I$.
\end{corollary}


The following proposition elucidates the structure of $P_I$ as $\W$
and $\heckeW{0}$-module.
\begin{proposition}
  \label{proposition.restriction.PI}
  Let $(-1)$ denote the sign representation of $\W$ as well as the
  corresponding representation of the Hecke algebra $\heckeW{0}$
  (sending $T_s$ to $-1$, or equivalently $\pi_s$ to $0$).
  \begin{itemize}
  \item[(a)] As a $\W$-module, $P_I \approx (-1)\uparrow_{\W_{I}}^{\W} $.
  \item[(b)] As an $\heckeW{0}$-module, $P_I \approx
    (-1)\uparrow_{\heckeW[\W_I]{0}}^{\heckeW{0}}$; it is a projective module.
  \item[(c)] The $P_I$'s are non isomorphic as $\heckeW{0}$-modules
    and thus as $\heckeWW$-modules.
  \end{itemize}
\end{proposition}

\begin{proof}
  First, for $s$ in $\W_I$, one has $v_I\act s = -v_I$. Moreover since
  ${}^I\W$ is a set of representative of left cosets for the quotient
  $\W_I\backslash\W$, the first point is proved.

  The second point is proved analogously, thanks to the fact that
  $v_I$ interpreted in the Hecke algebra as $\Phi(v_I)$ is an
  idempotent:
  \begin{equation}
    \Phi(v_I)^2 = 
    \left(\sum_{\nu \in \W_I} (-1)^{\lon(\nu)} \pi_\nu\right)^2 =
    \sum_{\nu \in \W_I} (-1)^{\lon(\nu)} \pi_\nu = \Phi(v_I)\,,
  \end{equation}
  since $\sum_{\nu \in \W_I} (-1)^{\lon(\nu)} \nu \pi_s = 0$ for any
  $s\in\W_I$. 

  The $P_I$'s being projective $\heckeW{0}$-modules, it is sufficient
  to prove that the associated semi-simple modules (obtained by
  factoring out by their radicals) are pairwise non-isomorphic.
  Namely, consider for each $I$
  \begin{equation}
    M_I := (-1)\uparrow_{H_I}^{H}/\rad\left((-1)\uparrow_{H_I}^{H}\right)\,,
  \end{equation}
  writing for short $H := \heckeW{0}$ and $H_I := \heckeW[\W_I]{0}$.
  Since $(-1)$ is simple as $H_I$-module,
  \begin{equation}
    M_I = (-1)\uparrow_{H_I/\rad(H_I)}^{H/\rad(H)}\, .
  \end{equation}
  Recall~\cite{Norton.1979} that $H/\rad(H)$ (resp. $H_I/\rad(H_I)$) is the
  commutative algebra over $(T_s)_{s\in S}$ (resp. $(T_s)_{s\in I}$).
  Hence, a simple $H/\rad(H)$-module $S^0_J$ is characterized by the set
  $J\in S$ of $s$ such that $\pi_s(S^0_J) = 0$. Since there is a simple
  $H$-module for each subset,
  \begin{equation}
    M_I \approx \bigoplus_{J\supset I} S^0_J\, .
  \end{equation}
  In particular, the $M_I$'s are pairwise non isomorphic $H$-modules,
  as desired.
\end{proof}

As a consequence of the preceding proof, if we denote $P^0_I$ the
projective $\heckeW{0}$-module associated with $S^0_I$, then, as an
$\heckeW{0}$-module
\begin{equation}
  P_I\approx \bigoplus_{J\supset I} P^0_J\,.
\end{equation}

In the following corollary, we focus on the particular type $A$, where a
$\sg{n}$-module (resp.  $H_n(0)$-module) is characterized by its so-called
characteristic, which is a symmetric function (resp.  a non-commutative
symmetric function).  Here is the characteristic of the modules $P_I$:
\begin{corollary}
  In type $A$, the characteristic of the $\sg{n}$-module $P_I$ is the
  symmetric function $e_K := e_{k_1}\cdots e_{k_l}$, where
  $K:=(k_1,\dots,k_l)$ is the composition associated to the complement
  of $I$.

  The characteristics of the $H_n(0)$-module $P_I$ is the
  noncommutative symmetric function $\Lambda^K:=\Lambda_{k_1}\cdots
  \Lambda_{k_l}$, where $K:=(k_1,\dots,k_l)$ is the composition
  associated to the complement of $I$.
\end{corollary}

\subsubsection{$\heckeWW$ as an incidence algebra}

Since $\{v_\sigma\}_{\sigma\in\W}$ is a basis of $\kW$, we may consider
for each $(\sigma,\tau)$ in $\W^2$ the matrix element
$\melem_{\s,\tau} \in \End(\kW)$ which maps $v_{\s'}$ to $v_\tau$ if
$\s'=\s$ and to $0$ otherwise.

We can now prove the main result of this section.
\begin{theorem}
  For each $(\sigma,\tau)$ in $\W^2$, let $\melem_{\s,\tau} \in \End(\kW)$
  defined by 
  \begin{equation}
    \melem_{\s,\tau}(\s') = \delta_{\s', \s}\, v_\tau\,. 
  \end{equation}
  Then the family $\{\melem_{\sigma,\tau} \suchthat \Rec(\sigma) \supset
  \Rec(\tau)\}$ is a vector space basis of $\heckeWW$.
\end{theorem}
\begin{proof}
  This family is free by construction. It has the appropriate size
  because for any finite Coxeter group
  $|\{\sigma \suchthat \Rec(\sigma)=I\}| =
  |\{\sigma \suchthat \Rec(\sigma)=S\backslash I\}|$
  (a bijection between the two sets is given by $u\mapsto \omega u$,
  where $\omega$ is the maximal element of $\W$; see~\cite[Exercise~10
  p.~57]{Bjorner_Brenti.2005}).
  
  It remains to check that any of its element $\melem_{\sigma,\tau}$
  preserves $i$-left antisymmetries and therefore is indeed in
  $\heckeWW$. Take $i\in S$, and consider the basis of $P_{\{i\}}$ of
  proposition~\ref{proposition.base.PI}: $\{ v_\sigma' \suchthat
  i\notin \Rec(\sigma') \}$; an element of this basis is either killed
  by $\melem_{\sigma,\tau}$ or sent to another element of this basis.
  Therefore, $P_{\{i\}}$ is stable by $\melem_{\sigma,\tau}$.
\end{proof}

Recall (see \cite{Stanley.1999.EnumerativeCombinatorics1}) that the
\emph{incidence algebra} $\k[P]$ of a poset $(P, \preccurlyeq)$ is the
algebra whose basis elements $e_{u,v}$ are indexed by the couples $(u,
v)\in P^2$ such $u\preccurlyeq v$, and whose multiplication rule is
given by:
\begin{equation}
  \label{eq.incidence.algebra.mult}
  e_{u,v} \cdot e_{u',v'} = \delta_{v,u'}\, e_{u,v'}\,.
\end{equation}
Here we need a slightly more general extension of this notion where $P$ is not
a partially ordered set but only a pre-order (not necessarily anti-symmetric). 

\begin{corollary}
  \label{corollary.incidenceAlgebra}
  $\heckeWW$ is isomorphic to the incidence algebra of the pre-order $(\W,
  \preccurlyeq)$ where
  \begin{equation}
    \s\preccurlyeq\s' 
    \quad\text{whenever}\quad 
    \Rec(\s)\supset \Rec(\s')\,.
  \end{equation}
\end{corollary}
\begin{proof}
  By construction, the $\melem_{\sigma,\tau}$ satisfy the usual
  product rule of incidence algebras:
  \begin{equation}
    \melem_{\sigma,\tau} \melem_{\sigma',\tau'} =
    \delta_{\tau,\sigma'} \melem_{\sigma,\tau'}.
  \end{equation}
  This is sufficient since $(\melem_{\sigma,\tau})$ is a basis. 
\end{proof}
The representation theory of $\heckeWW$ (projective and simple modules
and the Cartan matrix) will follow straightforwardly from this
corollary. Furthermore, a good way to think about the pre-order
$\preccurlyeq$ is to view it as the transitive closure of a graph $G$
which essentially encodes the action of the generators of $\heckeWW$.
Namely, define $G$ as the graph with vertex set $\W$ and where
$\s\rightarrow\s'$ is an edge whenever there exists $s\in S$ such that
$\s'=\s s$ and $\Rec(\s)\subset \Rec(\s')$ (cf.
Figure~\ref{figure.vsigma}).

\begin{figure}
  \begin{tikzpicture}[baseline=(current bounding box.east)]
    \tikzstyle{vertex}=[rectangle]
    \node(123) at (2,3) [vertex] {$v_{123} = 123 - 213 - 132 + 231 + 312 - 321$};
    \node(213) at (0,2) [vertex] {$v_{2\underline{1}3} = {213} - {312}$};
    \node(132) at (4,2) [vertex] {$v_{13\underline{2}} = {132} - {231}$};
    
    \node(231) at (0,1) [vertex] {$v_{23\underline{1}} = {231} - {321}$};
    \node(312) at (4,1) [vertex] {$v_{31\underline{2}} = {312} - {321}$};
    
    \node(321) at (2,0) [vertex] {$v_{3\underline{21}} = {321}$};

    \draw[<-](123) -- node [above left] {1} (213);
    \draw[<-](123) -- node [above right] {2} (132);
    \draw[<->](213) -- node [left] {2} (231);
    \draw[<->](132) -- node [right] {1} (312);
    \draw[<-](231) -- node [below left] {1} (321);
    \draw[<-](312) -- node [below right] {2} (321);
  \end{tikzpicture}

  \caption{The basis $(v_\sigma)_\sigma$ of $\ksg 3$, together with
    the graph structure which makes $\heckesg 3$ into an incidence
    algebra. Underlined: the recoils of the permutations.
  }
  \label{figure.vsigma}
\end{figure}

\begin{corollary}
  \label{theorem.projective}
  \begin{itemize}
  \item[(a)] The ideal $\melem_\s\heckeWW$ is isomorphic to
    $P_{S\backslash\Rec(\s)}$ as a right $\heckeWW$-module;
  \item[(b)] The idempotents $\melem_\s:=\melem_{\s,\s}$ give a
    maximal decomposition of the identity into orthogonal idempotents
    in $\heckeWW$;
  \end{itemize}
\end{corollary}



Hence we get the full description of the projective $\heckeWW$-modules.
\begin{corollary}
  The family of modules $(P_I)_{I\subset S}$ forms a complete set of
  representatives of the indecomposable projective modules of $\heckeWW$.
\end{corollary}

\begin{proof}
  The statement follows from (a) and (b) and Proposition~\ref{proposition.restriction.PI} (c).
\end{proof}

\subsubsection{Cartan matrix and the boolean lattice}

We now turn to the description of the Cartan matrix. For any $I\subset
S$, let $\alpha(I)$ be the shortest permutation such that
$\Rec(\alpha) = S\backslash I$ (the choice of the \emph{shortest} is
in fact irrelevant). For any $(I,J)$ such that $I \subset J \subset
S$, define $\melem_{I,J} := \melem_{\alpha(I),\alpha(J)}$.
  \begin{equation}
    \melem_{I,J} (v_\sigma) =
    \begin{cases}
      v_J & \text{if $\sigma=\alpha(I)$,} \\
      0 & \text{otherwise.}
    \end{cases}
  \end{equation}

\begin{remark}
  Using the incidence algebra structure, the sandwich $\melem_I
  \heckeWW \melem_J$ is $\k.\melem_{I,J}$ if $I\subset J$ and is
  trivial otherwise. 
\end{remark}

By Corollary~\ref{theorem.projective}, $\Hom(P_I, P_J)$ is isomorphic to the
sandwich $\melem_I \heckeWW \melem_J$. Therefore we have the following
corollary:
\begin{corollary}
  \label{proposition.heckesg.cartan}
  \begin{equation}
    \dim \Hom(P_I, P_J) = \dim \melem_I \heckeWW \melem_J=
    \begin{cases}
      1 & \text{if $I\supset J$,} \\
      0 & \text{otherwise.}
    \end{cases}
  \end{equation}
\end{corollary}
In other words, the Cartan matrix of $\heckeWW$ is the incidence
matrix of the boolean lattice. This suggests the existence of a close
relation between $\heckeWW$ and the incidence algebra of the boolean
lattice.

Recall that an algebra is called \emph{elementary} (or sometimes
\emph{reduced}) if its simple modules are all one dimensional.
Starting from an algebra $A$, it is possible to get a canonical
elementary algebra by the following process. Start with a maximal
decomposition of the identity $1=\sum_i e_i$ into orthogonal
idempotents.  Two idempotents $e_i$ and $e_j$ are \emph{conjugate} if
$e_i$ can be written as $a e_j b$ where $a$ and $b$ belongs to $A$, or
equivalently, if the projective modules $e_i A$ and $e_j A$ are
isomorphic. Select an idempotent $e_c$ in each conjugacy classes $c$
and put $e := \sum e_c$. Then, it is well
known~\cite{Curtis_Reiner.1990} that the algebra $eAe$ is elementary
and that the functor $M \mapsto M e$ which sends a right $A$-module to
a $eAe$-module is an equivalence of category. Recall finally that two
algebras $A$ and $B$ such that the categories of $A$-modules and
$B$-modules are equivalent are said \emph{Morita equivalent}.  Thus
$A$ and $eAe$ are Morita-equivalent. 
\begin{corollary}
  Let $e$ be the idempotent defined by $e := \sum_{I \subset S}
  \melem_I$. Then the algebra $e \heckeWW e$ is isomorphic to the
  incidence algebra $\k[B_S]$ of the boolean lattice $B_S$ of subsets
  of $S$. Consequently, $\heckeWW$ and $\k[B_S]$ are Morita
  equivalent.
\end{corollary}
Actually, the previous construction is fairly general when $A$ is the
incidence algebra of a pre-order $\preccurlyeq$, the construction of $eAe$
boils down to taking the incidence algebra of the canonical order $\leq$
associated to $\preccurlyeq$ by contracting its strongly connected components,
or equivalently by picking any representative. For $\heckeWW$, the strongly
connected components of the pre-order $(\W,\preccurlyeq)$ are by construction
$(\{\sigma \suchthat \Rec(\sigma)=I\})_{I\subset S}$, and the associated order
is simply the boolean lattice $B_S$ of subsets $I$ of $S$.


\subsubsection{Simple modules}

The simple modules are obtained as quotients of the projective modules
by their radical: $S_I:=P_I / \sum_{J\subsetneq I} P_J$. This has the
natural effect of making a vector $v$ in $S_I$ both $i$-antisymmetric
for $i\in I$ and $i$-symmetric for $i\not\in I$.
\begin{theorem}
  \label{theorem.heckesg.simple}
  The modules $(S_I)_{I\subset S}$ form a complete set of
  representatives of the simple modules of $\heckeWW$.  Moreover, the
  projection of the family $\{v_\s \suchthat \Rec(\s)=I\}$ in $S_I$ forms a
  vector space basis of $S_I$.
\end{theorem}

\begin{proof}
  This is a consequence of Proposition~\ref{proposition.heckesg.cartan}.
  Indeed, since there is no non trivial morphism from $P_I$ to itself, the
  radical of $P_I$ is the sum of the images of all morphisms from the $P_J$
  for $J \neq I$ to $P_J$. But there exists up to constant at most one such
  morphisms, that is when $I \supset J$, then $P_I\subset P_J$.  Hence
  $\rad(P_I) = \sum_{J\subsetneq I} P_J$. As a consequence, together with
  Proposition \ref{proposition.base.PI}, we get the given basis.
\end{proof}

The following proposition elucidates the structure of $S_I$ as $\W$
and $\heckeW{0}$-module.
\begin{proposition}
  \label{proposition.restriction.SI}
  Let $K$ be the composition associated to \FIXME{the complement of}$I$.
  \begin{itemize}
  \item[(a)] As a $\W$-module, $S_I$ is isomorphic to the subspace of
    $\kW$ of vectors which are simultaneously $i$-antisymmetric for
    $i$ in $I$ and $i$-symmetric for $i$ in $S\backslash I$. 
    In type $A$ this is the usual Young's
    representation $V_K$ indexed by the ribbon $K$, and each
    $v_\sigma$, with $\Des(\sigma)=I$ corresponds to the basis element
    of $V_K$ indexed by the standard ribbon tableaux associated to
    $\sigma^{-1}$. Its character is the ribbon Schur symmetric
    function $s_K$.
  \item[(b)] As a $\heckeW{0}$-module, $S_I$ is the indecomposable
    projective module indexed by $S\backslash I$. In type $A$, its
    character is the noncommutative symmetric function $R_K$.
  \end{itemize}
\end{proposition}
\begin{proof}
  (a) This description of $S_I$ is clear by construction. In type $A$
  this is a well known description of the ribbon representation.

  (b) For $\s$ such that $\Rec(\s) = S\backslash I$, define $b_\s$
  as the image of $v_I \act \s$ by the canonical quotient map.  Then,
  $\{b_\s \suchthat \Rec(\s) = S\backslash I\}$ is a basis of
  $S_I$. Furthermore, one can easily check that:
  \begin{equation}
    b_\s \pi_s =
    \begin{cases}
      b_\s     & \text{if $s\in \Des(s)$,}\\
      b_{\s s} & \text{if $s\not\in \Des(s)$ and $\Rec(\s s)=S\backslash I$,}\\
      0        & \text{otherwise.}
    \end{cases}
  \end{equation}
  One recognizes the usual description of the projective
  $\heckeW{0}$-module $P_{S\backslash I}$ of
  Section~\ref{section.reptheo0Hecke}.
\end{proof}
Note that $S_I$ is also isomorphic to the simple module $S_{S\backslash I}$ of
the transpose algebra $\heckeWW^*$.

\subsubsection{Induction, restriction, and Grothendieck rings}

In this subsection, we concentrate on the tower of type-$A$ algebras
$(\heckesg n)_n$.

Let $\calG := \calG\left((\heckesg{n})_n\right)$ and $\calK:=
\calK\left((\heckesg{n})_n\right)$ be respectively the Grothendieck
rings of the characters of the simple and projective modules of the
tower of algebras $(\heckesg{n})_n$.  Let furthermore $\cartan$ be the
Cartan map from $\calK$ to $\calG$. It is the algebra and coalgebra
morphism which gives the projection of a module onto the direct sum of
its composition factors. It is given by
\begin{equation}
  \cartan(P_I) = \sum_{I\supset J} S_J\,.
\end{equation}
Since the indecomposable projective modules are indexed by
compositions, it comes out as no surprise that the structure of
algebras and coalgebras of $\calG$ and $\calK$ are isomorphic to
$\qsym$ and $\ncsf$, respectively.  However, we do not get Hopf
algebras, because the structures of algebras and coalgebras are not
compatible.

\begin{proposition}
  The following diagram gives a complete description of the structures
  of algebras and of coalgebras on $\calG$ and $\calK$:
  \begin{equation}
    \entrymodifiers={+<20pt>}
    \vcenter{
    \xymatrix@R=0pt@C=1.1cm{
      (\qsym, .)
      & ({\calG},.)
      \ar@{^(->>}[l]_{\chi(S_I)\mapsto M_{I}}
      & ({\calK},.)
      \ar@{^(->>}[l]_{\cartan}
      \ar@{^(->>}[r]^{\chi(P_I)\mapsto F_{I}}
      & (\qsym, .)
      \\
      (\ncsf, \Delta)
      & ({\calG},\Delta) 
      \ar@{^(->>}[l]^{\chi(S_I)\mapsto R_{I}}
      & ({\calK},\Delta) 
      \ar@{^(->>}[l]^{\cartan}
      \ar@{^(->>}[r]_{\chi(P_I)\mapsto \Lambda^{I^c}}
      & (\ncsf, \Delta)
    }}
  \end{equation}
\end{proposition}
\begin{proof}
  The bottom line is already known from Proposition
  \ref{proposition.restriction.PI} and the fact that, for all $m$ and
  $n$, the following diagram commutes:
\begin{equation}
  \entrymodifiers={+<10pt>}
  \vcenter{\xymatrix{
      {\hecke{m}{0} \otimes \hecke{n}{0}\ } \ar@{^(->}[d] \ar@{^(->}[r] &
      {\ \hecke{m+n}{0}} \ar@{^(->}[d]
      \\
      {\heckesg{m} \otimes \heckesg{n}\ } \ar@{^(->}[r] &
      {\ \heckesg{m+n}}
    }}
\end{equation}
Thus the map which sends a module to the characteristic of its
restriction to $\hecke{n}{0}$ is a coalgebra morphism. The isomorphism
from $({\calK},.)$ to $\qsym$ is then obtained by Frobenius duality
between induction of projective modules and restriction of simple
modules. And the last case is obtained by applying the Cartan map $\cartan$.
\end{proof}
It is important to note that the algebra $({\calG},.)$ is not the dual
of the coalgebra $(\calK, \Delta)$ because the dual of the restriction
of projective modules is the so-called \emph{co-induction} of simple
modules which, in general, is not the same as the induction for non
self-injective algebras.

Finally, we briefly describe the Grothendieck rings for the adjoint
algebra $\heckeWW^*$ which preserves symmetries. The projective
modules are defined by the two following equivalent formulas:
\begin{equation}
  \label{eq.def.P.I.star}
  P_I := \bigcap_{s \in I} \ker (1-\ls)
  = \left(\sum_{\nu \in \W_I} \nu\right)\act \heckeWW^*\ .
\end{equation}
As $\W$-module (resp. $\heckeW{0}$-module), they are isomorphic to
the modules induced by the trivial modules of $\W_I$ (resp.
$\heckeW[\W_I]{0}$) whose Frobenius characteristic are complete
symmetric functions. The rest of our arguments can be adapted easily,
yielding the following diagram:
  \begin{equation}
    \entrymodifiers={+<20pt>}
    \vcenter{
    \xymatrix@R=0pt@C=1.1cm{
      (\qsym, .)
      & ({\calG},.)
      \ar@{^(->>}[l]_{\chi(S_I)\mapsto X_{I}}
      & ({\calK},.)
      \ar@{^(->>}[l]_{\cartan}
      \ar@{^(->>}[r]^{\chi(P_I)\mapsto F_{I^c}}
      & (\qsym, .)
      \\
      (\ncsf, \Delta)
      & ({\calG},\Delta) 
      \ar@{^(->>}[l]^{\chi(S_I)\mapsto R_{I^c}}
      & ({\calK},\Delta) 
      \ar@{^(->>}[l]^{\cartan}
      \ar@{^(->>}[r]_{\chi(P_I)\mapsto S^{I^c}}
      & (\ncsf, \Delta)
    }},
  \end{equation}
where $(X_I)_I$ is the dual basis of the elementary basis $(\Lambda_I)_I$ of
$\ncsf$. Thus we have a representation theoretical interpretation of many
bases of $\ncsf$ and $\qsym$.


\section{The algebra of non-decreasing functions}
\label{section.ndf}

\begin{definition}
  Let $\ndf{n}$ be the set of \emph{non-decreasing functions} from
  $\{1,\dots,n\}$ to itself. Its cardinality is given by:
  \begin{equation}
    \label{equation.card.ndf}
    \binom{2n-1}{n-1} = \sum_{k=1}^n \binom n k \binom{n-1}{k-1}\,,
  \end{equation}
  where, on the right, functions are counted according to the size $k$
  of their images.

  The composition and the neutral element $\id_n$ make $\ndf{n}$ into
  a monoid and we denote by $\ndfa{n}$ its monoid algebra.  The monoid
  $\ndf{n}\times \ndf{m}$ can be identified as the submonoid of
  $\ndf{n+m}$ whose elements stabilize both $\{1,\dots,n\}$ and
  $\{n+1,\dots,n+m\}$. This makes $(\ndfa{n})_n$ into a tower of
  algebras.
\end{definition}


\TODO{Recontacter Pr Putcha; cf mail 23 mars}

The semi-group properties of $\ndf{n}$ have been studied (see
e.g.~\cite{Gomes_Howie.1992}, where $\mathcal O_n$ coincides with
$\ndf{n}$ striped of the identity). In
particular~\cite[Theorem~4.8]{Gomes_Howie.1992}, one can take as
idempotent generators for $\ndf{n}$ the functions $\pi_i$ et $\opi_i$
defined by:
\begin{equation}
  \begin{aligned}
    \pi_i(i+1):=i  &\quad\text{ and } \quad \pi_i(j):=j, \text{ for } j\ne i+1,\\
    \opi_i(i):=i+1 &\quad\text{ and } \quad \pi_i(j):=j, \text{ for } j\ne i.
  \end{aligned}
\end{equation}
The functions $\pi_i$ satisfy the braid relations, together with a new
relation:
\begin{equation}
  \label{equ.pres.ndpf}
  \pi_i^2=\pi_i
  \qquad \text{and} \qquad
  \pi_{i+1}\pi_i\pi_{i+1} = \pi_i\pi_{i+1}\pi_i = \pi_{i+1}\pi_i\, .
\end{equation}

This readily defines a morphism $\phi: \pi_{\hecke{n}{0}} \mapsto \pi_{\ndfa
  n}$ of $\hecke{n}{0}$ into $\ndfa n$. Its image is the monoid algebra of
\emph{non-decreasing parking functions} which will be discussed in
Section~\ref{section.ndpf} and of which Equation~(\ref{equ.pres.ndpf}) actually
gives a presentation. The same properties hold for the operators $\opi_i$'s.
Although this is not a priori obvious, it will turn out that the two morphisms
$\phi:\pi_{\hecke{n}{0}}\mapsto \pi_{\ndfa n}$ and $\overline \phi:
\opi_{\hecke{n}{0}}\mapsto\opi_{\ndfa n}$ are compatible, making $\ndfa n$
into a quotient of $\heckesg{n}$
(Proposition~\ref{proposition.ndfa.quotient}). This will be used in
Subsection~\ref{section.ndf.representationTheory} to deduce the representation
theory of $\ndfa{n}$.

\subsection{Representation on exterior powers, and link with $\heckesg n$}

We now want to construct a suitable faithful representation of
$\ndfa{n}$ where the existence of the epimorphism from $\heckesg{n}$
onto $\ndfa{n}$ becomes clear.

The \emph{natural representation} of $\ndfa{n}$ is obtained by taking the
vector space $\k^n$ with canonical basis $e_1,\dots,e_n$, and letting a
function $f$ act on it by $e_i.f=e_{f(i)}$. For $n>2$, this representation is
a faithful representation of the monoid $\ndf{n}$ but not of the algebra, as
$\dim \ndfa{n}=\binom{2n-1}{n-1} \gg n^2$. However, since $\ndf{n}$ is a
monoid, the diagonal action on \emph{exterior powers}
\begin{equation}
  \label{eq.action.exterieur}
  (x_1\wedge\dots\wedge x_k) \act f :=
  (x_1 \act f)\wedge\dots\wedge (x_k \act f)
\end{equation}
still defines an action. Taking the \emph{exterior powers} $\bigwedge^k
\k^n$ of the natural representation gives a new representation, whose
basis $\{e_S:=e_{s_1}\wedge\dots\wedge e_{s_k}\}$ is indexed by
subsets $S=\{s_1,\dots,s_k\}$ of $\{1,\dots,n\}$. The action of a
function $f$ in $\ndf{n}$ is simply given by (note the absence of
sign!):
\begin{equation}
  e_S.f =
  \begin{cases}
    e_{f(S)} & \text{ if $|f(S)|=|S|$},\\
    0      & \text{ otherwise.}
  \end{cases}
\end{equation}

We call \emph{representation of $\ndfa{n}$ on exterior powers} the
representation of $\ndfa{n}$ on $\bigoplus_{k=1}^n \bigwedge^k \k^n$,
which is of dimension $2^n-1$ (it turns out that we do not need to
include the component $\bigwedge^0 \k^n$ for our purposes).
\begin{lemma}
  \label{lemma.ndfa.faithfull}
  For $n>0$, the representation 
  $\bigoplus_{k=1}^n \bigwedge^k \k^n$ of $\ndfa{n}$ is faithful.
\end{lemma}
\begin{proof}
  We exhibit a triangularity property. For a function $f$ in
  $\ndf{n}$, let $\im f:=f(\{1,\dots,n\})$ be its image set, and
  $R(f)$ be the preimage of $\im f$ defined by $R(f):=\{ \min \{x \suchthat
  f(x) = y\},\ y \in S \}$. Put a partial order on $\ndfa{n}$ by
  setting $f\prec g$ if $\im f=\im g$ and $R(f)$ is lexicographically
  smaller than $R(g)$. Fix a function $f$. If the representation
  matrix of a function $g$ has coefficient $1$ on row $\im f$ and
  column $R(f)$, that is if $g(R(f))=\im f$, then $g\preceq f$.
  
  Remark: for $n>0$, the component $\bigwedge^0 \k^n$ is not needed
  because $\im f$ and $R(f)$ are non empty.
\end{proof}

We now realize the representation of $\ndfa{n}$ on the $k$-th exterior
power as a representation of $\heckesg{n}$. Let us start from the
$k$-th exterior product $\bigwedge^k\k^n$ of the natural
representation of $\sg{n}$ on $\k^n$. It is isomorphic to
$V_{(n-k,1,\dots,1)} \bigoplus V_{(n-k+1,1,\dots,1)}$, where
$V_\lambda$ denotes the irreducible representation of $\sg n$ indexed
by the partition $\lambda$. It can be realized as a submodule
of the regular representation of $\sg{n}$ using the classical Young
construction by mean of the row-symmetrizers and
column-antisymmetrizers on the skew ribbon
\begin{equation}
  {\def\lr#1{\multicolumn{1}{|c|}{#1}}
    \begin{array}{ccccc}
      \cline{1-1}
      \lr{k}\\
      \cline{1-1}
      \lr{\vdots}\\
      \cline{1-1}
      \lr{1}\\
      \cline{1-4}
      &\lr{k+1}&\lr{\cdots}&\lr{n}\\
      \cline{2-4}
    \end{array}
  }\,.
\end{equation}
This only involves left antisymmetries on the values $1,\dots,k-1$ and
symmetries on the values $k+1,\dots,n-1$. Therefore, $\bigwedge^k\k^n$
can alternatively be realized as the $\heckesg n$-module
\begin{equation}
  P_n^k := P_{ \{1,\dots,k-1\} }\quad /\quad \sum_{s=k+1}^{n-1} P_{\{1,\dots,k-1,s\}}
\end{equation}
(an element of $P_{\{1,\dots,k-1\}}$ is by construction left
antisymmetric on the values $1,\dots,k-1$, and the quotient by each
$P_{\{1,\dots,k-1,s\}}$ makes it $s$-left symmetric). In general, this
construction turns any $\sg n$-module indexed by a shape made of
disconnected rows and columns into an $\heckesg n$-module; it does not
apply to shapes like $\lambda=(2,1)$, as they involve symmetries or
antisymmetries between non-consecutive values, like $1$ and $3$.

A basis of $P_n^k$ indexed by subsets of size $k$ of $\{1,\dots,n\}$
is obtained by taking for each such subset $S$ the image in the
quotient $P_n^k$ of
\begin{equation}
  e_S := \sum_{
  \substack{
    \sigma(S) = \{1,\dots,k\}\\
    \sigma(i) < \sigma(j), \text{ for $i,j \notin S$ with $i<j$}
  }}
  (-1)^{\lon (\sigma)} \sigma\, .
\end{equation}
It is straightforward to check that the actions of $\pi_i$ and
$\opi_i$ of $\heckesg{n}$ on $e_S$ of $P_n^k$ coincide with the actions
of $\pi_i$ and $\opi_i$ of $\ndfa n$ on $e_S$ of $\bigwedge^k \k^n$
(justifying a posteriori the identical notations). In the sequel, we
identify the modules $P_n^k$ and $\bigwedge^k \k^n$ of $\heckesg{n}$
and $\ndfa n$, and we call \emph{representation on exterior powers of
  $\heckesg{n}$} its representation on $\bigoplus_{k=1}^n \bigwedge^k
\k^n$.

Using Lemma~\ref{lemma.ndfa.faithfull} we are in position to state the
following
\begin{proposition}
  \label{proposition.ndfa.quotient}
  $\ndfa{n}$ is the quotient of $\heckesg{n}$ obtained by considering
  its representation on exterior powers. The restriction of this
  representation of $\heckesg{n}$ to $\sga{n}$, $\hecke{n}{0}$, and
  $\hecke{n}{-1}$ yield respectively the usual representation of
  $\sg{n}$ on exterior powers, the algebra $\ndpf n$ of non-decreasing
  parking functions (see Section~\ref{section.ndpf}), and the
  Temperley-Lieb algebra $\TL{n}{-1}$.
\end{proposition}

\begin{proof}
  \def\tl{\operatorname{e}}%
  The only case which remains is $q=-1$. Recall
  that~\cite{Temperley_Lieb.1971} the Temperley-Lieb algebra
  $\TL{n}{q}$ is the quotient of the Hecke algebra by the relations
  \begin{equation}
    \label{eq.def.temperley-lieb}
    \tl_i\tl_{i\pm1}\tl_i = q \tl_i\,,
  \end{equation}
  where $\tl_i=T_i+q$. As $\ndpf n$, its dimension is the Catalan number
  $C_n$.

  The algebra $A$ is generated by the operators $\tl_i = \pi_i-\opi_i$, which
  satisfy the relations:
  \begin{equation}
    \begin{gathered}
      \tl_i^2 = 0\,,\\
      \tl_i \tl_{i\pm1} \tl_i = -\tl_i\,.\\
    \end{gathered}
  \end{equation}
  Therefore $A$ is a quotient of the Temperley-Lieb algebra~$\TL{n}{-1}$.
  We now prove that the quotient is trivial by exhibiting $C_n$
  elements of $A$ satisfying a triangularity property w.r.t. those of
  $\ndpf n$.

  Let $\leq$ be the pointwise partial order on $\ndf n$ such that
  $f\leq g$ if and only if $f(i)\leq g(i)$ for all $i$. The following
  properties are easily verified:
  \begin{equation}
    \begin{gathered}
      f  \pi_i \leq f \leq f \opi_i\\
      f\leq g \quad\Longrightarrow\quad f \pi_i \leq g \pi_i\\
      f\leq g \quad\Longrightarrow\quad f \opi_i \leq g \opi_i
    \end{gathered}
  \end{equation}
  Take $i_1,\dots,i_k$ such that the product $\pi_{i_1}\cdots\pi_{i_k}$
  is reduced in $\ndpf n$, and consider a function
  \begin{equation}
    f=\pi_{i_1}\cdots\pi_{i_{k_1}}\opi_{i_{k_1+1}}\cdots\opi_{i_{k_2}}\pi_{i_{k_2+1}}\cdots
  \end{equation}
  of $\ndf n$ appearing in the expansion of the product
  \begin{equation}
    \tl_{i_1}\cdots\tl_{i_k}=(\pi_{i_1}-\opi_{i_1})\cdots(\pi_{i_1}-\opi_{i_1})\,.
  \end{equation}
  Clearly, $\pi_{i_1}\cdots\pi_{i_k} \leq f$. Furthermore, if
  equality holds then
  \begin{equation}
    \pi_{i_1}\cdots\pi_{i_k}\leq\pi_{i_1}\cdots\pi_{i_{k_1}}\pi_{i_{k_2+1}}\cdots
    \leq f = \pi_{i_1}\cdots\pi_{i_k}\,.
  \end{equation}
  Since the product is reduced, $k_1=k$, and therefore,
  $\pi_{i_1}\cdots\pi_{i_k}$ appears with coefficient $1$ in
  $\tl_{i_1}\cdots\tl_{i_k}$.
\end{proof}

\subsection{Representation theory}
\label{section.ndf.representationTheory}

In this section, we derive the representation theory of $\ndf n$ from
that of $\heckesg n$. An alternative more combinatorial approach would
be to construct by inclusion/exclusion a graded basis $( \gr f )_{f\in
  \ndf n}$ of $\ndfa f$ such that:
\begin{equation}
  \gr f \gr g =
  \begin{cases}
    \gr {fg} & \text{if $|\im f|=|\im g|=|\im fg|$,}\\
    0 & \text{otherwise.}
  \end{cases}
\end{equation}
Then, looking at the principal modules $\gr f.\ndfa n=\k.\{\gr
g\suchthat \fibers(g)=\fibers(f)\}$ splits the regular representation
into a direct summand of $\binom{n-1}{k-1}$ copies of $P_n^k$ for each
$k$.  Combinatorial proofs for the Cartan matrix and the simple
modules are then lengthy but straightforward.

\subsubsection{Projective modules, simple modules, and Cartan  matrix}

Let $\delta$ be the usual homology border map ($\delta^2=0$):
\begin{equation}
  \delta:
  \begin{cases}
    P_n^k & \to P_n^{k-1} \\
    S:=\{s_1,\dots,s_k\} & \mapsto \sum_{i\in \{1,\dots,k\}}
    (-1)^{k-i} S \backslash \{s_i\}
  \end{cases}\, ,
\end{equation}
which induces the following exact sequence of morphisms of
$\ndfa n$-modules:
\begin{equation}\label{equation.long.exact.pi}
  \def\arrow{\stackrel\delta\longrightarrow}
  0 \rightarrow P_n^n \arrow \cdots \arrow
  P_n^{k+1} \arrow P_n^{k} \arrow P_n^{k-1} \arrow
  \cdots \arrow P_n^1 \arrow P_n^0 \rightarrow 0\,/
\end{equation}
%
%
%
%
For $k=1,\dots, n$, set $S_n^k := P_n^k / \ker \delta$, so that we
have the short exact sequence:
\begin{equation}
  0 \rightarrow S_n^{k+1} \rightarrow P_n^k \rightarrow S_n^k
  \rightarrow 0\, .
\end{equation}
In particular, $\dim S_n^k=\binom{n-1}{k-1}$ since $\dim P_n^k =\binom n k$.
The following proposition states that the morphism of
Equation~(\ref{equation.long.exact.pi}) are essentially the only non trivial
morphisms between the $(P_n^k)_{k=1,\dots,n}$.
\begin{proposition}
  \label{proposition.ndfa.hom.pi}
  Let $k$ and $l$ be two integers in $\{1,\dots,n\}$. Then,
  \begin{equation}
    \dim \Hom(P_n^k, P_n^l) =
    \begin{cases}
      1 & \text{ if } k\in \{l, l+1\} \,, \\
      0 & \text{ otherwise.}
    \end{cases}
  \end{equation}  
\end{proposition}
 
\begin{proof}
  From the Cartan matrix of $\heckesg n$ (see
  Proposition~\ref{proposition.heckesg.cartan}) we can deduce the dimension of
  $\Hom(P_n^k, P_n^l)$. Indeed, any non trivial $\ndfa n$-morphism $\phi$ from
  $P_n^k$ to $P_n^l$ is an $\heckesg n$-morphism, and thus can be lifted to a
  non trivial $\heckesg n$-morphism $\psi$ from $P_{\{1,\dots,k-1\}}$ to the
  projective module $P_{\{1,\dots,l-1\}}$.

  If $k<l$, $\dim \Hom(P_{\{1,\dots,k-1\}}, P_{\{1,\dots,l-1\}})= 0$,
  and therefore $\dim \Hom(P_n^k, P_n^l)=0$. Otherwise, $\dim
  \Hom(P_{\{1,\dots,k-1\}}, P_{\{1,\dots,l-1\}})= 1$.  If $k> l+1$,
  then $P_{\{1,\dots,k-1\}}$ is mapped by $\psi$ to its unique copy in
  $P_{\{1,\dots,l, l+1\}}$, and is therefore killed in the quotient
  \begin{equation}
    P_n^l := P_{ \{1,\dots,l-1\} }\quad /\quad \sum_{s=l+1}^{n-1} P_{\{1,\dots,l-1,s\}}\,.
  \end{equation}
  therefore $\dim \Hom(P_n^k, P_n^l)=0$. In the remaining cases we can
  conclude that $\Hom(P_n^k, P_n^k)=\k.\id$, and $\Hom(P_n^k,
  P_n^{k-1})=\k.\delta$.
\end{proof}

We are now in position to describe the projective and simple modules.
\begin{proposition}
  \label{proposition.ndfa.representationTheory}
  The modules $(P_n^k)_{k=1,\dots,n}$ form a complete set of
  representatives of the indecomposable projective modules of $\ndfa
  n$.

  The modules $(S_n^k)_{k=1,\dots,n}$ form a complete set of
  representatives of the simple modules of $\ndfa n$.
\end{proposition}

\begin{proof}
  It follows from Proposition~\ref{proposition.ndfa.hom.pi} that the
  modules $(P_n^k)_{k=1,\dots,n}$ are both indecomposable ($\dim
  \Hom(P_n^k, P_n^k)=1$) and non isomorphic ($\dim \Hom(P_n^k,
  P_n^l)=0$ if $k<l$). It remains to prove (i) that each of them is
  projective and (ii) that we obtain all the projective modules this
  way. It then follows from the description of the morphisms between
  the $P_n^k$ that the $(S_n^k)_{k=1,\dots,n}$ form a complete set of
  representatives of the simple modules.

  We first achieve (i) by constructing explicitly an idempotent
  $e_n^k$ such that the principal ideal $e_n^k \ndfa n$ is isomorphic
  to $P_n^k$. Define the idempotent $e_n^k$ as follows:
  \begin{equation}
    \label{equation.def.indemp.ndfa}
    e_n^k := \pi_{n-1}\pi_{n-2}\cdots\pi_{k+1}\pi_{k}\
       (1-\pi_{k-1})(1-\pi_{k-2})\cdots(1-\pi_{2})(1-\pi_{1})\,.
  \end{equation}
  To prove that $e_n^k$ is indeed an idempotent, it is sufficient to
  use the presentation of $\ndfa n$ given in
  Equation~(\ref{equ.pres.ndpf}) to check that
  \begin{equation}
    \label{equation.idemp.idemp}
    e_n^k\pi_i =
      \begin{cases}
        e_n^k & \text{if $i>=k$,}\\
        0     & \text{otherwise.}
  \end{cases}
  \end{equation}
  Alternatively, we could have shown that $e_n^k$ is the image in
  $\ndfa n$ of the classical hook idempotent of the $0$-Hecke algebra
  \begin{equation}
    \prod_j \pi_j  \prod_j (1-\pi_j)\,,
  \end{equation}
  where the left (resp. on the right) product ranges over a reduced
  word of the maximal permutation of the parabolic subgroup
  $\sg{k,\dots,n-1}$ (resp. $\sg{1,\dots k-1}$).

  We want now to prove that $P_n^k$ is isomorphic to $e_n^k \ndfa n$.
  Let us denote $B_n^k$ the set of functions in $\ndf{n}$ which are
  injective on $\{1,\dots,k\}$ and such that $f(i)=f(k)$ for all
  $i\geq k$.  It is clear that for any function $f$ of $\ndf n$, the
  actual value of $e_n^k f$ depends only on the set $f(\{1,\dots,
  k\})$. Moreover, for any $f\in B_n^k$,
  \begin{equation}
    \label{equation.triangle.idemp}
    e_n^kf = f + \sum \text{functions strictly smaller than $f$}\,,
  \end{equation}
  where the order considered on the set of functions is the product order on
  the tuple $(f(1),\dots,f(n))$ (ie. $f\geq g$ iff $f(i) \geq g(i)$ for all
  $i$). It follows by triangularity that $\{e_n^k f\ \suchthat\ f\in B_n^k\}$ is a
  basis for $e_n^k \ndfa n$.

  Now the linear map
  \begin{equation}
    \phi_n^k:
    \begin{cases}
      e_n^k \ndfa n & \to P_n^{k} \\
      e_n^k f & \mapsto \{f(1), \dots, f(k)\} 
    \end{cases}\, ,
  \end{equation}
  is in fact a morphism of $\ndfa n$-module. Indeed,
  \begin{equation}
    \phi_n^k(e_n^k f \pi_i) = \phi_n^k(e_n^k (f \pi_i)) 
       = \{f(1)\pi_i, \dots, f(k)\pi_i\} = \{f(1), \dots, f(k)\} \pi_i\,,
  \end{equation}
  and the same holds for $\opi_i$. Since the cardinal of $B_n^k$ is the same
  as the dimension of $P_n^k$, on gets that $\phi_n^k$ is in fact an
  isomorphism. Therefore $P_n^k$ is projective.

  To derive (ii), note that the faithful $\ndfa n$-module
  $\bigoplus_{k=1}^n \bigwedge^k \k^n = \bigoplus_{k=1}^n P^k_n$
  (Lemma~\ref{lemma.ndfa.faithfull}) is now projective. And it is a
  general property that any projective module of a finite dimensional
  algebra $A$ occurs as a submodule of any faithful projective module
  $M$ of $A$ (the right regular representation of $A$ is a sub-module
  of a certain number of copy of $M$, and the decomposition of any of
  its projective module is known to be unique up to an isomorphism).




\end{proof}
Hence Proposition \ref{proposition.ndfa.hom.pi} actually gives the Cartan
matrix of $\ndfa n$.
\subsubsection{Induction, restriction, and Grothendieck rings}

\begin{proposition}
  The restriction and induction of indecomposable projective modules
  and simple modules are described by:
  \begin{equation}
    P_{n_1+n_2}^k \downarrow^{\ndfa{n_1+n_2}}_{\ndfa{n_1}\otimes \ndfa{n_2}} \approx
    \bigoplus_{\substack{
        n_1+n_2=n\\
        k_1+k_2=k\\
        1\leq k_i\leq n_i \text{ or } k_i=n_i=0
      }}
    P_{n_1}^{k_1} \otimes P_{n_2}^{k_2}
  \end{equation}
  \begin{equation}
    P_{n_1}^{k_1} \otimes P_{n_2}^{k_2} \uparrow^{\ndfa{n_1+n_2}}_{\ndfa{n_1}\otimes \ndfa{n_2}}\approx
    P_{n_1+n_2}^{k_1+k_2} \oplus P_{n_1+n_2}^{k_1+k_2-1}
  \end{equation}
  \begin{equation}
    S_{n_1+n_2}^k \downarrow^{\ndfa{n_1+n_2}}_{\ndfa{n_1}\otimes \ndfa{n_2}} =
    \bigoplus_{
      \substack{
        n_1+n_2=n\\
        k_1+k_2 \in \{k,k+1\}\\
        1\leq k_i\leq n_i \text{ or } k_i=n_i=0
      }
    }
    S_{n_1}^{k_1} \otimes S_{n_2}^{k_2}
  \end{equation}
  \begin{equation}
    S_{n_1}^{k_1} \otimes S_{n_2}^{k_2} \uparrow^{\ndfa{n_1+n_2}}_{\ndfa{n_1}\otimes \ndfa{n_2}}\approx
    S_{n_1+n_2}^{k_1+k_2}
  \end{equation}
\end{proposition}

Those rules yield structures of commutative algebras and cocommutative
coalgebras on the Grothendieck rings of $\ndfa{n}$. However, we do not
get Hopf algebras, because the structures of algebras and coalgebras
are not compatible (the coefficient of $\chi(P_1^1)\otimes\chi(P_1^1)$
differs in $\Delta(\chi(P_1^1)\chi(P_1^1))$ and
$\Delta(\chi(P_1^1))\Delta(\chi(P_1^1))$).

\begin{proposition}
  The Grothendieck rings $\calG$ and $\calK$ of $\ndfa{n}$ can be
  realized as quotients or subcoalgebras of $\sym$, $\qsym$, and
  $\ncsf$, as described in the following diagram:
  \begin{equation}
    \vcenter{
    \entrymodifiers={+<12pt>}
    \xymatrix@R=10pt@C=0.6cm{
      (\sym, .)
      \ar@{->>}[rr]^{h_\lambda \mapsto \chi\left(S_{|\lambda|}^{\lon(\lambda)}\right)}
      && (\calG,.)
      & (\calK,.)
      \ar@{^(->>}[l]_{\cartan}
      &&& (\ncsf, .)
      \ar@{->>}[lll]_{R_I \mapsto \chi\left(P_{|I|}^{\lon(I)}\right)}
      \\
      (\qsym, \Delta)
      \ar@{->>}[rr]_{F_I \mapsto \chi\left(S_{|I|}^{\lon(I)}\right)}
      && (\calG,\Delta) \ar@{<-->}[ru]^{*}
      & (\calK,\Delta) 
      \ar@{^(->>}[l]^{\cartan}
      \ar@{^(->}[rrr]_{\chi(P_n^k)\mapsto { \sum_{\lambda\partof n,\lon(\lambda)=k} m_\lambda}}
      &&& (\sym, \Delta)
    }}
  \end{equation}
\end{proposition}

\section{The algebra of non-decreasing parking functions}
\label{section.ndpf}


\begin{definition}
  A \emph{nondecreasing parking function} of size $n$ is a
  nondecreasing function $f$ from $\{1,2,\dots n\}$ to $\{1,2,\dots
  n\}$ such that $f(i) \leq i$, for all $i\leq n$.
  
  The composition of maps and the neutral element $\id_n$ make the set
  of nondecreasing parking function of size $n$ into a monoid denoted
  $\ndpf{n}$.
\end{definition}
Parking functions where introduced
in~\cite{Konheim_Weiss.1966.ParkingFunctions}. It is well known that
the nondecreasing parking functions are counted by the Catalan numbers
$C_n = \frac1{n+1}\binom{2n}{n}$. It is also clear that $\ndpf{n}$ is
the sub-monoid of $\ndf{n}$ generated by the $\pi_i$'s.

\subsection{Representation theory}

\subsubsection{Simple modules}

The goal of the sequel is to study the representation theory of
$\ndpf{n}$, or equivalently of its algebra $\ndpfa{n}$. The following
remark allows us to deduce the representations of $\ndpfa{n}$ from the
representations of $\hecke{n}{0}$.
\begin{proposition}
  The kernel of the algebra epi-morphism $\phi : \hecke{n}{0} \to \ndpfa{n}$
  defined by $\phi(\pi_i) = \pi_i$ is a sub-ideal of the radical of
  $\hecke{n}{0}$.
\end{proposition}
\begin{proof}
  It is well known (see \cite{Norton.1979}) that the quotient of $\hecke{n}{0}$ by its
  radical is a commutative algebra. Consequently, $\pi_i\pi_{i+1}\pi_i -
  \pi_i\pi_{i+1} = [\pi_i\pi_{i+1}, \pi_i]$ belongs to the radical of
  $\hecke{n}{0}$.
\end{proof}
As a consequence, taking the quotient by their respective radical
shows that the projection $\phi$ is an isomorphism from
$\hecke{n}{0}/\rad(\hecke{n}{0})$ to
$\ndpfa{n}/\rad(\ndpfa{n})$. Moreover $\ndpfa{n}/\rad(\ndpfa{n})$ is
isomorphic to the commutative algebra generated by the $\pi_i$ such
that $\pi_i^2=\pi_i$. As a consequence, $\hecke{n}{0}$ and $\ndpfa{n}$
share, roughly speaking, the same simple modules:

\begin{corollary}
  There are $2^{n-1}$ simple $\ndpfa{n}$-modules $S_I$, and they are
  all one dimensional. The structure of the module $S_I$, generated by
  $\eta_I$, is given by
\begin{equation}
\left\{
\begin{array}{r@{\,}cl}
\eta_I\act\pi_i &= 0      & \text{if $i\in I$,}\\
\eta_I\act\pi_i &= \eta_I & \text{otherwise.}
\end{array}
\right.
\end{equation}
\end{corollary}

\subsubsection{Projective modules}

The projective modules of $\ndpf{n}$ can be deduced from the ones of
$\ndf{n}$.
\begin{theorem}
  Let $I= \{s_1, \dots, s_k\} \subset \{1,\dots,n-1\}$. Then, the principal
  sub-module
  \begin{equation}
    P_I :=
    (e_1 \wedge e_{s_1+1} \wedge \dots \wedge e_{s_k+1}) \act \ndpfa{n}
    \quad \subset\quad \bigwedge^{k+1} \k^n
  \end{equation}
  is an indecomposable projective module. Moreover, the set $(P_I)_{I
    \compof n}$ is a complete set of representatives of indecomposable
  projective modules of $\ndpfa{n}$.
\end{theorem}


This suggests an alternative description of the algebra $\ndpfa{n}$.
Let $G_{n,k}$ be the lattice of subsets of $\{1,\dots,n\}$ of size $k$
for the \emph{product order} defined as follows. Let $S :=
\{s_1<s_2<\dots<s_k\}$ and $T := \{t_1<t_2<\dots<t_k\}$ be two
subsets.  Then,
\begin{equation}
  S \leq_G T
  \qquad\text{if and only if}\qquad
  s_i \leq t_i \text{, for $i=1,\dots, k$.}
\end{equation}
One easily sees that $S \leq_G T$ if and only if there exists a
nondecreasing parking function $f$ such that $e_S = e_T \act f$. This
lattice appears as the Bruhat order associated to the Grassman
manifold $G^n_k$ of $k$-dimensional subspaces in $\CC^n$.
\begin{theorem}
  There is a natural algebra isomorphism
  \begin{equation}
    \ndpfa{n}\ \approx\ \bigoplus_{k=0}^{n-1} \k[G_{n-1, k}]\,.
  \end{equation}
\end{theorem}
In particular the Cartan map $\cartan : \calK \to \calG$ is given by
the lattice $\leq_G$:
\begin{equation}
  \cartan(P_I)\ =\ \sum_{J,\ \Des(J) \leq_G \Des(I)} S_J\,.
\end{equation}
On the other hand, due to the commutative diagram
\begin{equation}
  \entrymodifiers={+<10pt>}
  \vcenter{\xymatrix{
      {\hecke{m}{0} \otimes \hecke{n}{0}\ } \ar@{->>}[d] \ar@{^(->}[r] &
      {\ \hecke{m+n}{0}} \ar@{->>}[d]
      \\
      {\ndpf{m} \otimes \ndpf{n}\ } \ar@{^(->}[r] &
      {\ \ndpf{m+n}}
    }}
\end{equation}
it is clear that the restriction of simple modules and the induction
of indecomposable projective modules follow the same rule as for
$\hecke{n}{0}$. The induction of simple modules can be deduced via the
Cartan map, giving rise to a new basis $G_I$ of $\ncsf$. The restriction of 
indecomposable projective modules leads to a new operation on compositions,
which seems not to be related to anything previously known. 
All of this is summarized by the following diagram:
  \begin{equation}
    \entrymodifiers={+<12pt>}
    \vcenter{
    \xymatrix@R=10pt@C=0.77cm{
      (\ncsf, .)
      && ({\calG},.)
      \ar@{^(->>}[ll]_{\chi(S_I)\mapsto G_{I}}
      & ({\calK},.)
      \ar@{^(->>}[l]_{\cartan}
      \ar@{^(->>}[rr]^{\chi(P_I)\mapsto R_{I}}
      && (\ncsf, .)
      \\
      (\qsym, \Delta)
      && (\calG,\Delta)
      \ar@{^(->>}[ll]^{\chi(S_I)\mapsto F_I}
      & (\calK,\Delta)
      \ar@{^(->>}[l]^{\cartan}
      \ar@{^(->>}[rr]_{\chi(P_I)\mapsto ???}
      && ???
    }}
\end{equation}

\section{Alternative constructions of $\heckesg{n}$ in type A}
\label{section.alternativeDefinitions}

In type $A$, the actions of the operators $s_i$ and $\pi_i$ on
permutations of $\sg n$ extend straightforwardly to an action on the
set $A^n$ of words $w$ of length $n$ over any totally ordered alphabet
$A$ by
\begin{equation}
  w \act \pi_i =
  \begin{cases}
    w      & \text{if $w_i \geq w_{i+1}$,}\\
    w s_i & \text{otherwise}.
  \end{cases}
\end{equation}
We may again construct the algebra $\k[s_i,\pi_i]$, and wonder whether
it is strictly larger than $\heckesg n$.
\begin{theorem}
  Let $A$ be a totally ordered alphabet of size at least $n$. Then the
  subalgebra of $\End(\k.A^n)$ generated by both sets of operators
  $\{s_i, \pi_i\}_{i=1,\dots,n-1}$ is isomorphic to $\heckesg n$.
\end{theorem}

\newcommand{\EndEvaluation}{\End_e}

This theorem is best restated and proved using a commuting property
with the monoid of non decreasing functions.  Recall that the
evaluation $e(w)$ of a word $w$ on an alphabet $A$ is the function
which counts the number of occurrences in $w$ of each letter $a\in A$;
for example, the evaluation of a permutation is the constant function
$a\mapsto 1$. For a given evaluation $e$, write $C_e$ the subspace of
$\k.A^n$ spanned by the words with evaluation $e$, and $p_e$ the
orthogonal projection on $C_e$.
Let finally $\EndEvaluation(\k.A^n)=\bigoplus_e \End(C_e)$ denote the
algebra of \emph{evaluation preserving endomorphisms} of $\k.A^n$.

\begin{theorem}
  \label{tilt}
  Let $A$ be a totally ordered alphabet of size at least $n$.
  Consider the monoid $\ndf A$ of non decreasing functions from $A$ to
  $A$, acting on values on the words of $A^n$. Then the commutant of
  $\ndf A$ in $\EndEvaluation(\k.A^n)$ coincides with $\k[s_i,\pi_i]$
  and is isomorphic to $\heckesg n$.

  Alternatively, $\k[s_i,\pi_i]$ is the commutant of the subalgebra of
  $\End(\k.A^n)$ generated by both $\ndf A$ and $(p_e)_e$.
\end{theorem}
\begin{proof}
  First it is clear that both $s_i$ and $\pi_i$ preserve the
  evaluation and commute with non decreasing functions; so
  $\k[s_i,\pi_i]$ is a subset of the commutant.

  Taking $n$ distinct letters in $A$ which we may call $1<\dots<n$
  yields a component $C_{\std}$ isomorphic to $\ksg n$, which we call
  \emph{standard component}. The restriction of $\k[s_i,\pi_i]$ on
  $C_{\std}$ is of course $\heckesg n$.

  Fix now some operator $f\in \EndEvaluation(\k.A^n)$ which commutes
  with $\ndf A$.  For each $i$, take a function $\pi_i$ in $\ndf A$
  such that for $1\leq j < k\leq n$, $\pi_i(j)=\pi_i(k)$ if and only
  if $j=i$ and $k=i+1$. A vector $v$ in the standard component is in
  the kernel of $\pi_i$ if and only if $v$ is $i$-left antisymmetric.
  Therefore, the restriction on $C_{\std}$ of $f$ preserves
  $i$-antisymmetries, and thus coincides with some operator $g$ of
  $\heckesg n$. On the other hand, for any component $C_e$ there
  exists a non decreasing function $f_e$ which maps $C_{\std}$ onto $C_e$.
  Since $f$ commutes with $f_e$, the action of $f$ on $C_e$ is
  determined by its action of $C_{\std}$, that is by $g$.
\end{proof}
Note: it is in fact sufficient to consider just 
$p_\std$ instead of $(p_e)_e$ in theorem~\ref{tilt}. 

The argument can in fact be generalized to any subset $W$ of words
containing some words with all letters distinct, and stable
simultaneously by the right action of $S_n$ and by the left action of
some monoid of non decreasing functions large enough to contain
analogues of the $\pi_i$'s and $f_e$'s. The following proposition
gives two typical examples of that situation
(A function $f$ from $\{1,\dots,n\}$ to itself is \emph{initial} if there
exists $k\leq n$ such that $\im(f) = \{1,\dots, k\}$; for parking functions
see~\cite{Konheim_Weiss.1966.ParkingFunctions})

\begin{proposition}
  Let $A$ be the totally ordered alphabet $\{1<\dots<n\}$.
  \begin{itemize}
  \item[(a)] Consider the monoid $\ndpf n$ of non decreasing parking
    functions, acting on the left on the set $\pf n$ of parking
    functions from $\{1,\dots,n\}$ to $\{1,\dots,n\}$. Then the
    commutant of $\ndpf n$ in $\EndEvaluation(\k.\pf n)$ coincides
    with $\k[s_i,\pi_i]$ and is isomorphic to $\heckesg n$.
  \item[(a)] Consider the monoid $\ndinitial n$ of non decreasing
    initial functions, acting on the left on the set $\initial n$ of
    initial functions from $\{1,\dots,n\}$ to $\{1,\dots,n\}$. Then
    the commutant of $\ndinitial n$ in $\EndEvaluation(\k.\initial n)$
    coincides with $\k[s_i,\pi_i]$ and is isomorphic to $\heckesg n$.
  \end{itemize}
\end{proposition}

\section{Research in progress}
\label{section.researchInProgress}

A first direction of research concerns the links between Hecke group
algebras and affine Hecke algebras~\cite{Hivert_Schilling_Thiery.HeckeGroupAffine.2007,Hivert_Schilling_Thiery.HeckeGroupAffine.2008}. It turns out that for any Weyl
group, $\heckeWW$ is the natural quotient of the generic untwisted affine Hecke
algebra $\affineheckeW{q}$, through its level zero action. On one
hand, this may shed some new light on the representation theory of the
affine Hecke algebras. On the other hand, this yields several
representations of $\heckesg n$ on the ring of polynomials
$\k[x_1,\dots,x_n]$, the $\pi_i$ acting, for example, by elementary sorting
on monomials or as isobaric divided differences. Furthermore,
the maximal commutative subring $\k[Y_1,\dots,Y_n]$ of
$\affinehecke{q}$, specializes to a maximal commutative subring of
dimension $n!$ of $\heckesg{n}$ by evaluating the symmetric functions
$\sym(Y_1,\dots,Y_n)$ on the alphabet $\frac{1-q^n}{1-q}$.  This is
likely to give some explicit description of the orthogonal idempotents
$\melem_\s$ of Theorem~\ref{theorem.projective}. More importantly,
this opens the door for links with the non symmetric Macdonald
polynomials (which can, among other things be defined as eigenvalues
of the $Y_i$~\cite{Kirillov_Noumi.1998.AffineHecke}).

The properties of the Hecke group algebras, and in particular their
representation theory seem also to generalize nicely to infinite
Coxeter groups, up to some little adaptations. First, when a parabolic
subgroup $W_I$ is infinite, it is not possible to realize the
projective module $P_I$ inside $\Wa$, or at least not without an
appropriate completion (because of the infinite alternating sum).
Reciprocally, $P_{S\backslash I}$ is not distinct anymore from
$\bigcup_{J\supset I} P_J$ (there are no element of $\W$ with descent
set $I$; cf also proposition~\ref{proposition.base.PI}).  This
suggests that $\heckeWW$ is Morita equivalent to the poset algebra of
some convex subset of the boolean lattice.

A last direction of research is the generalization of
sections~\ref{section.ndf} and~\ref{section.ndpf}. This essentially
boils down to the following question: what is the natural definition
of the representation on exterior powers for a general Coxeter group?
One such attempt in type $B$ gives rise to some tower of
self-injective monoid of signed non decreasing parking functions whose
sizes appear to be given by sequence A086618 of the encyclopedia of
integer sequences~\cite{Sloane}.

\bibliographystyle{alpha}
\bibliography{main}
\end{document}